\documentclass[amssymb,12pt]{article}
\newcount\localcountcount
\newcount\localboxcount
\def\newlocalcount#1{
	\advance\localcountcount by -1
	\allocationnumber= \localcountcount
	\ifnum\count10<\allocationnumber
	  \else
	    \errmessage{No room for a new local count.}
	  \fi
	\countdef#1=\allocationnumber
    }

\def\newlocalbox#1{
	\advance\localboxcount by -1
	\allocationnumber= \localboxcount
	\ifnum\count14<\allocationnumber
	  \else
	    \errmessage{No room for a new local box.}
	  \fi
	\chardef#1=\allocationnumber
    }

\newcount\coordinatescalex
\newcount\coordinatescaley
\def\bigger(#1/#2){
	\multiply\coordinatescalex by #1
	\divide\coordinatescalex by #2
	\multiply\coordinatescaley by #1
	\divide\coordinatescaley by #2
    }
\def\wider(#1/#2){
	\multiply\coordinatescalex by #1
	\divide\coordinatescalex by #2
    }
\def\taller(#1/#2){
	\multiply\coordinatescaley by #1
	\divide\coordinatescaley by #2
    }

\newlength{\objectmargin}
\newlength{\arrowmargin}
\newlength{\cornersize}
\newlength{\arrowheadlength}
\newlength{\nudgesize}
\def\enlarge #1 by #2/#3;{
	\multiply\csname #1\endcsname by #2
	\divide\csname #1\endcsname by #3
    }

\newif\ifflippinglabels

\newif\ifepimorphic

\newif\ifmonomorphic

\newcount\labelposition
\newcount\nudgex
\newcount\nudgey
\def\nudge(#1,#2){
	\global\advance\nudgex by #1
	\global\advance\nudgey by #2
    }
\def\nudgeright{
	\global\multiply\nudgex by \unitlength
	\global\advance\nudgex by \nudgesize
	\global\divide\nudgex by \unitlength
    }
\def\nudgeleft{
	\global\multiply\nudgex by \unitlength
	\global\advance\nudgex by -\nudgesize
	\global\divide\nudgex by \unitlength
    }
\def\nudgeup{
	\global\multiply\nudgey by \unitlength
	\global\advance\nudgey by \nudgesize
	\global\divide\nudgey by \unitlength
    }
\def\nudgedown{
	\global\multiply\nudgey by \unitlength
	\global\advance\nudgey by -\nudgesize
	\global\divide\nudgey by \unitlength
    }

\def\negate(#1){
	\csname #1\endcsname= -\csname #1\endcsname
    }

\def\makepositive(#1){
	\ifnum\csname #1\endcsname<0
	    \negate(#1)
	  \fi
    }

\def\calcgdivisor(#1,#2) into #3 {{
	\newlocalcount\left
	\newlocalcount\right
	\newlocalcount\temp

	\left=  #1
	\right= #2
	\makepositive(left)
	\makepositive(right)

	\loop
	    \ifnum\left<\right
		\temp=  \left
		\left=  \right
		\right= \temp
	      \fi
	  \ifnum\right>0 
	    \advance\left by -\right
	  \repeat
	\global\csname #3\endcsname= \left
    }}

\def\newpoint #1 {
	\expandafter\newlocalcount\csname #1x\endcsname
	\expandafter\newlocalcount\csname #1y\endcsname
    }

\def\setpoint #1 to (#2,#3){
	\csname #1x\endcsname = #2
	\csname #1y\endcsname = #3
    }

\def\copypoint #1 to #2 {
	\csname #2x\endcsname = \csname #1x\endcsname
	\csname #2y\endcsname = \csname #1y\endcsname
    }

\def\addpoint #1 to #2 {
	\advance\csname #2x\endcsname by \csname #1x\endcsname
	\advance\csname #2y\endcsname by \csname #1y\endcsname
    }

\def\subtractpoint #1 from #2 {
	\advance\csname #2x\endcsname by -\csname #1x\endcsname
	\advance\csname #2y\endcsname by -\csname #1y\endcsname
    }

\def\multiplypoint #1 by #2 {
	\multiply\csname #1x\endcsname by #2
	\multiply\csname #1y\endcsname by #2
    }

\def\dividepoint #1 by #2 {
	\divide\csname #1x\endcsname by #2
	\divide\csname #1y\endcsname by #2
    }

\def\printpoint #1 {
     \message{(\number\csname #1x\endcsname,\number\csname #1y\endcsname)}
    }

\def\lookuptrig #1:#2 {
	\ifnum1=#1
	    \ifnum0=#2		% 1:0
		\multiply\trigpointx by 10000
		\multiply\trigpointy by 0
	      \else		% 1:1
		\multiply\trigpointx by 7071
		\multiply\trigpointy by 7071
	      \fi
	  \fi
	\ifnum2=#1		% 2:1
	    \multiply\trigpointx by 8944
	    \multiply\trigpointy by 4472
	  \fi	
	\ifnum3=#1
	    \ifnum1=#2		% 3:1
		\multiply\trigpointx by 9487
		\multiply\trigpointy by 3162
	      \else		% 3:2
		\multiply\trigpointx by 8321
		\multiply\trigpointy by 5547
	      \fi
	  \fi	
	\ifnum4=#1
	    \ifnum1=#2		% 4:1
		\multiply\trigpointx by 9701
		\multiply\trigpointy by 2425
	      \else		% 4:3
		\multiply\trigpointx by 8000
		\multiply\trigpointy by 6000
	      \fi
	  \fi	
    }

\def\calcunsignedtrig #1:#2 {
	\ifnum#1<#2
	    \lookuptrig #2:#1
	    \advance\trigpointx by -\trigpointy
	    \advance\trigpointy by \trigpointx
	    \advance\trigpointx by -\trigpointy
	    \trigpointx= -\trigpointx
	  \else
	    \lookuptrig #1:#2
	  \fi
    }

\def\calctrig #1:#2 {
	\ifnum0>#2
	    \ifnum0>#1
		\calcunsignedtrig -#1:-#2
	        \trigpointx= -\trigpointx
	      \else
		\calcunsignedtrig #1:-#2
	      \fi
	    \trigpointy= -\trigpointy
	  \else
	    \ifnum0>#1
		\calcunsignedtrig -#1:#2
	        \trigpointx= -\trigpointx
	      \else
		\calcunsignedtrig #1:#2
	      \fi
	  \fi
    }

\def\calctrigoffset #1 along #2:#3 {
	\trigpointx= #1
	\divide\trigpointx by \unitlength
	\trigpointy= \trigpointx
	\calctrig #2:#3
	\dividepoint trigpoint by 10000
    }

\def\trigoffset #1 by #2 along #3:#4 {
	\calctrigoffset #2 along #3:#4
	\addpoint trigpoint to #1
    }

\def\newvector #1 {
	\newpoint {#1head}
	\newpoint {#1tail}
	\newpoint {#1slope}
	\expandafter\newlocalcount\csname #1length\endcsname
    }

\def\shortenhead #1 by #2 {
	\subtractpoint {#2} from {#1head}
	\newlocalcount\changeinlength
	\ifnum\csname #1slopex\endcsname=0
	    \changeinlength= \csname #2y\endcsname
	  \else
	    \changeinlength= \csname #2x\endcsname
	  \fi
	\makepositive(changeinlength)
	\advance\csname #1length\endcsname by -\changeinlength
    }

\def\shortentail #1 by #2 {
	\addpoint {#2} to {#1tail}
	\newlocalcount\changeinlength
	\ifnum\csname #1slopex\endcsname=0
	    \changeinlength= \csname #2y\endcsname
	  \else
	    \changeinlength= \csname #2x\endcsname
	  \fi
	\makepositive(changeinlength)
	\advance\csname #1length\endcsname by -\changeinlength
    }

\def\setvector #1 from #2 to #3 {
	\copypoint #2 to {#1tail}
	\copypoint #3 to {#1head}

	\copypoint #3 to {#1slope}
	\subtractpoint #2 from {#1slope}

	\ifnum\csname #1slopex\endcsname=0
	    \csname #1length\endcsname= \csname #1slopey\endcsname
	  \else
	    \csname #1length\endcsname= \csname #1slopex\endcsname
	  \fi
	\makepositive(#1length)

	\newlocalcount\gdivisor
	\calcgdivisor(
		\csname #1slopex\endcsname,
		\csname #1slopey\endcsname
	    ) into gdivisor
	\dividepoint {#1slope} by {\gdivisor}
    }

\def\getvectormidpoint #1 into #2
    {
	\copypoint {#1head} to #2
	\subtractpoint  {#1tail} from #2
	\multiplypoint #2 by {\labelposition}
	\dividepoint #2 by 100
	\addpoint {#1tail} to #2
	\global\labelposition=50
    }

\def\drawline #1 {
	\put(\csname #1tailx\endcsname,\csname #1taily\endcsname){\line( \csname #1slopex\endcsname,\csname #1slopey\endcsname){\csname #1length\endcsname}}
    }

\def\drawsecondhead #1 {{
	\newpoint {}
	\newpoint trigpoint
	\copypoint {#1head} to {}
	\trigoffset {}
		by {-\arrowheadlength}
		along {\csname #1slopex\endcsname}:{\csname #1slopey\endcsname}
\put(\x,\y){\vector(\csname #1slopex\endcsname,\csname #1slopey\endcsname){0}}
    }}

\def\drawarrowtail #1 {{
	\newpoint {}
	\newpoint trigpoint
	\copypoint {#1tail} to {}
	\trigoffset {}
		by {\arrowheadlength}
		along {\csname #1slopex\endcsname}:{\csname #1slopey\endcsname}
\put(\x,\y){\vector(\csname #1slopex\endcsname,\csname #1slopey\endcsname){0}}
    }}

\def\drawvector #1 {
	\put(\csname #1tailx\endcsname,\csname #1taily\endcsname){\vector( \csname #1slopex\endcsname,\csname #1slopey\endcsname){\csname #1length\endcsname}}
	\ifepimorphic
	    \drawsecondhead {#1}
	    \global\epimorphicfalse
	  \fi	
	\ifmonomorphic
	    \drawarrowtail {#1}
	    \global\monomorphicfalse
	  \fi	
    }

\def\drawdoublevector #1 {
	\drawvector {#1}
	\put(\csname #1tailx\endcsname,\csname #1taily\endcsname){\vector(-\csname #1slopex\endcsname,-\csname #1slopey\endcsname){ 0 }}
    }

\def\expandxbounds(#1){
	\ifnum\leastobjectx>#1
	    \global\leastobjectx=#1
	  \fi
	\ifnum\greatestobjectx<#1
	    \global\greatestobjectx=#1
	  \fi
    }

\def\expandybounds(#1,#2){
	\ifnum\leastobjecty>#2
	    \global\leastobjecty=#2
	  \fi
	\ifnum\greatestobjecty<#1
	    \global\greatestobjecty=#1
	  \fi
    }

\def\object #1 #2,#3 #4;{
	\expandafter\newlocalcount\csname objectx#1\endcsname
	\expandafter\newlocalcount\csname objecty#1\endcsname
	\csname objectx#1\endcsname=#2
	\csname objecty#1\endcsname=#3
	\multiply\csname objectx#1\endcsname by \coordinatescalex
	\multiply\csname objecty#1\endcsname by \coordinatescaley

	\expandafter\newlocalcount\csname halfwidth#1\endcsname
	\expandafter\newlocalcount\csname halfheight#1\endcsname
	{
	    \newlocalbox\objectbox

	    \savebox{\objectbox}{$\displaystyle{#4}$}
	    \ifdim\ht\objectbox=0pt 
		\savebox{\objectbox}{$\bullet$}
	      \fi

	    \newpoint {}
	    \setpoint {} to (#2,#3)
	    \multiply\x by \coordinatescalex
	    \multiply\y by \coordinatescaley

	    \put(\x,\y){\makebox(0,0){\usebox{\objectbox}}}
	
	    \newlocalcount\halfwidth
	    \newlocalcount\halfheight

	    \expandxbounds(\x)

	    \halfwidth=  \wd\objectbox
	    \halfheight= \ht\objectbox
	    \advance\halfheight by \dp\objectbox

	    \divide\halfwidth  by 2
	    \divide\halfheight by 2

	    \divide\halfwidth  by \unitlength
	    \divide\halfheight by \unitlength

	    \newlocalcount\top
	    \newlocalcount\bottom
	    \top= \y
	    \bottom= \y
	    \advance\top    by \halfheight
	    \advance\bottom by -\halfheight
	    \expandybounds(\top,\bottom)

	    \newlocalcount\margin
	    \margin= \objectmargin
	    \divide\margin by \unitlength
	    \advance\halfwidth  by \margin
	    \advance\halfheight by \margin

	    \global\csname halfwidth#1\endcsname=  \halfwidth
	    \global\csname halfheight#1\endcsname= \halfheight
	}
    }

\newlength{\starsize}
\def\starzero #1;{{
	\newpoint center
	\newpoint foo
	\newpoint trigpoint
	\centerx= \csname objectx#1\endcsname
	\centery= \csname objecty#1\endcsname
	\copypoint center to foo
	\trigoffset foo by {\starsize} along 0:1
	\put(\foox,\fooy){\makebox(0,0){$\bullet$}}
	\copypoint center to foo
	\trigoffset foo by {\starsize} along 0:-1
	\put(\foox,\fooy){\makebox(0,0){$\bullet$}}
	\copypoint center to foo
	\trigoffset foo by {\starsize} along 1:0
	\put(\foox,\fooy){\makebox(0,0){$\bullet$}}
	\copypoint center to foo
	\trigoffset foo by {\starsize} along -1:0
	\put(\foox,\fooy){\makebox(0,0){$\bullet$}}
    }}

\def\starone #1;{{
	\newpoint center
	\newpoint foo
	\newpoint trigpoint
	\centerx= \csname objectx#1\endcsname
	\centery= \csname objecty#1\endcsname
	\copypoint center to foo
	\trigoffset foo by {\starsize} along 1:1
	\put(\foox,\fooy){\makebox(0,0){$\bullet$}}
	\copypoint center to foo
	\trigoffset foo by {\starsize} along 1:-1
	\put(\foox,\fooy){\makebox(0,0){$\bullet$}}
	\copypoint center to foo
	\trigoffset foo by {\starsize} along -1:1
	\put(\foox,\fooy){\makebox(0,0){$\bullet$}}
	\copypoint center to foo
	\trigoffset foo by {\starsize} along -1:-1
	\put(\foox,\fooy){\makebox(0,0){$\bullet$}}
    }}

\def\startwo #1;{{
	\newpoint center
	\newpoint foo
	\newpoint trigpoint
	\centerx= \csname objectx#1\endcsname
	\centery= \csname objecty#1\endcsname
	\copypoint center to foo
	\trigoffset foo by {\starsize} along 1:2
	\put(\foox,\fooy){\makebox(0,0){$\bullet$}}
	\copypoint center to foo
	\trigoffset foo by {\starsize} along 1:-2
	\put(\foox,\fooy){\makebox(0,0){$\bullet$}}
	\copypoint center to foo
	\trigoffset foo by {\starsize} along -1:2
	\put(\foox,\fooy){\makebox(0,0){$\bullet$}}
	\copypoint center to foo
	\trigoffset foo by {\starsize} along -1:-2
	\put(\foox,\fooy){\makebox(0,0){$\bullet$}}
	\copypoint center to foo
	\trigoffset foo by {\starsize} along 2:1
	\put(\foox,\fooy){\makebox(0,0){$\bullet$}}
	\copypoint center to foo
	\trigoffset foo by {\starsize} along -2:1
	\put(\foox,\fooy){\makebox(0,0){$\bullet$}}
	\copypoint center to foo
	\trigoffset foo by {\starsize} along 2:-1
	\put(\foox,\fooy){\makebox(0,0){$\bullet$}}
	\copypoint center to foo
	\trigoffset foo by {\starsize} along -2:-1
	\put(\foox,\fooy){\makebox(0,0){$\bullet$}}
    }}

\def\starthree #1;{{
	\newpoint center
	\newpoint foo
	\newpoint trigpoint
	\centerx= \csname objectx#1\endcsname
	\centery= \csname objecty#1\endcsname
	\copypoint center to foo
	\trigoffset foo by {\starsize} along 1:3
	\put(\foox,\fooy){\makebox(0,0){$\bullet$}}
	\copypoint center to foo
	\trigoffset foo by {\starsize} along 1:-3
	\put(\foox,\fooy){\makebox(0,0){$\bullet$}}
	\copypoint center to foo
	\trigoffset foo by {\starsize} along -1:3
	\put(\foox,\fooy){\makebox(0,0){$\bullet$}}
	\copypoint center to foo
	\trigoffset foo by {\starsize} along -1:-3
	\put(\foox,\fooy){\makebox(0,0){$\bullet$}}
	\copypoint center to foo
	\trigoffset foo by {\starsize} along 3:1
	\put(\foox,\fooy){\makebox(0,0){$\bullet$}}
	\copypoint center to foo
	\trigoffset foo by {\starsize} along -3:1
	\put(\foox,\fooy){\makebox(0,0){$\bullet$}}
	\copypoint center to foo
	\trigoffset foo by {\starsize} along 3:-1
	\put(\foox,\fooy){\makebox(0,0){$\bullet$}}
	\copypoint center to foo
	\trigoffset foo by {\starsize} along -3:-1
	\put(\foox,\fooy){\makebox(0,0){$\bullet$}}
	\copypoint center to foo
	\trigoffset foo by {\starsize} along 2:3
	\put(\foox,\fooy){\makebox(0,0){$\bullet$}}
	\copypoint center to foo
	\trigoffset foo by {\starsize} along 2:-3
	\put(\foox,\fooy){\makebox(0,0){$\bullet$}}
	\copypoint center to foo
	\trigoffset foo by {\starsize} along -2:3
	\put(\foox,\fooy){\makebox(0,0){$\bullet$}}
	\copypoint center to foo
	\trigoffset foo by {\starsize} along -2:-3
	\put(\foox,\fooy){\makebox(0,0){$\bullet$}}
	\copypoint center to foo
	\trigoffset foo by {\starsize} along 3:2
	\put(\foox,\fooy){\makebox(0,0){$\bullet$}}
	\copypoint center to foo
	\trigoffset foo by {\starsize} along -3:2
	\put(\foox,\fooy){\makebox(0,0){$\bullet$}}
	\copypoint center to foo
	\trigoffset foo by {\starsize} along 3:-2
	\put(\foox,\fooy){\makebox(0,0){$\bullet$}}
	\copypoint center to foo
	\trigoffset foo by {\starsize} along -3:-2
	\put(\foox,\fooy){\makebox(0,0){$\bullet$}}
    }}

\def\starfour #1;{{
	\newpoint center
	\newpoint foo
	\newpoint trigpoint
	\centerx= \csname objectx#1\endcsname
	\centery= \csname objecty#1\endcsname
	\copypoint center to foo
	\trigoffset foo by {\starsize} along 1:4
	\put(\foox,\fooy){\makebox(0,0){$\bullet$}}
	\copypoint center to foo
	\trigoffset foo by {\starsize} along 1:-4
	\put(\foox,\fooy){\makebox(0,0){$\bullet$}}
	\copypoint center to foo
	\trigoffset foo by {\starsize} along -1:4
	\put(\foox,\fooy){\makebox(0,0){$\bullet$}}
	\copypoint center to foo
	\trigoffset foo by {\starsize} along -1:-4
	\put(\foox,\fooy){\makebox(0,0){$\bullet$}}
	\copypoint center to foo
	\trigoffset foo by {\starsize} along 4:1
	\put(\foox,\fooy){\makebox(0,0){$\bullet$}}
	\copypoint center to foo
	\trigoffset foo by {\starsize} along -4:1
	\put(\foox,\fooy){\makebox(0,0){$\bullet$}}
	\copypoint center to foo
	\trigoffset foo by {\starsize} along 4:-1
	\put(\foox,\fooy){\makebox(0,0){$\bullet$}}
	\copypoint center to foo
	\trigoffset foo by {\starsize} along -4:-1
	\put(\foox,\fooy){\makebox(0,0){$\bullet$}}
	\copypoint center to foo
	\trigoffset foo by {\starsize} along 3:4
	\put(\foox,\fooy){\makebox(0,0){$\bullet$}}
	\copypoint center to foo
	\trigoffset foo by {\starsize} along 3:-4
	\put(\foox,\fooy){\makebox(0,0){$\bullet$}}
	\copypoint center to foo
	\trigoffset foo by {\starsize} along -3:4
	\put(\foox,\fooy){\makebox(0,0){$\bullet$}}
	\copypoint center to foo
	\trigoffset foo by {\starsize} along -3:-4
	\put(\foox,\fooy){\makebox(0,0){$\bullet$}}
	\copypoint center to foo
	\trigoffset foo by {\starsize} along 4:3
	\put(\foox,\fooy){\makebox(0,0){$\bullet$}}
	\copypoint center to foo
	\trigoffset foo by {\starsize} along -4:3
	\put(\foox,\fooy){\makebox(0,0){$\bullet$}}
	\copypoint center to foo
	\trigoffset foo by {\starsize} along 4:-3
	\put(\foox,\fooy){\makebox(0,0){$\bullet$}}
	\copypoint center to foo
	\trigoffset foo by {\starsize} along -4:-3
	\put(\foox,\fooy){\makebox(0,0){$\bullet$}}
    }}

\def\getobjectposition #1 into #2 {
	\csname #2x\endcsname = \csname objectx#1\endcsname
	\csname #2y\endcsname = \csname objecty#1\endcsname
    }

\def\getobjectsize #1 into #2 {
	\csname #2x\endcsname = \csname halfwidth#1\endcsname
	\csname #2y\endcsname = \csname halfheight#1\endcsname
    }

\def\labelbox(#1,#2)(#3,#4)#5{{
	\newlocalbox\labelbox
	\savebox{\labelbox}{${#5}$}

	\newpoint {}
	\setpoint {} to (#1,#2)

	\newlocalcount\margin
	\margin= \arrowmargin
	\divide\margin by \unitlength
	
	\newlocalcount\height
	\newlocalcount\width
	\newlocalcount\depth
	\height= \ht\labelbox
	\width=  \wd\labelbox
	\depth=  \dp\labelbox
	\divide\height by \unitlength
	\divide\width  by \unitlength
	\divide\depth  by \unitlength

	\ifnum#3=0
	    \ifflippinglabels
	    	\advance\x by \margin
	      \else
	    	\advance\x by -\width
	    	\advance\x by -\margin
	      \fi
	     \advance\y by \y
	     \advance\y by -\height
	     \advance\y by \depth
	     \divide\y by 2
	  \else
	    \newlocalcount\slopesign
	    \slopesign= #3
	    \multiply\slopesign by #4

	    \ifflippinglabels
		\advance\y by -\height
		\advance\y by -\margin
		\multiply\slopesign by -1
	      \else
		\advance\y by \depth
		\advance\y by \margin
	      \fi
	    
	    \ifnum\slopesign>0
		\advance\x by -\width
		\advance\x by -\margin
	      \else
		\ifnum\slopesign<0
		    \advance\x by \margin
		  \else
		    \divide\width by 2
		    \advance\x by -\width
		  \fi
	      \fi
	  \fi

	\put(\x,\y){\usebox{\labelbox}}

	\newlocalcount\top
	\newlocalcount\bottom
	\top=    \y
	\bottom= \y
	\advance\top    by \height
	\advance\bottom by -\depth
	\expandybounds(\top,\bottom)
	\global\flippinglabelsfalse
    }}

\def\avoidobject #1 (#2,#3){{
	\getobjectsize {#1} into delta

	\newlocalcount\hweight
	\newlocalcount\vweight

	\hweight=\deltay
	\vweight=\deltax
	\multiply\hweight by #2
	\multiply\vweight by #3
	\makepositive(hweight)
	\makepositive(vweight)

	\ifnum\hweight<\vweight
		\ifnum#3<0 \deltay= -\deltay \fi
		\deltax=\deltay
		\multiply\deltax by #2
		\divide\deltax by #3
	    \else
		\ifnum#2<0 \deltax= -\deltax \fi
		\deltay=\deltax
		\multiply\deltay by #3
		\divide\deltay by #2
	    \fi

	\global\deltax= \deltax
	\global\deltay= \deltay
    }}

\def\calcarrow #1 #2 {
	\newpoint start
	\newpoint end
	\getobjectposition {#1} into start
	\getobjectposition {#2} into end

	\newvector arrow
	\setvector arrow from start to end

	\newpoint delta
	\avoidobject #1 (\arrowslopex,\arrowslopey)
	\shortentail arrow by delta

	\avoidobject #2 (\arrowslopex,\arrowslopey)
	\shortenhead arrow by delta

	\addpoint nudge to arrowhead
	\addpoint nudge to arrowtail
	
	\global\nudgex=0
	\global\nudgey=0
    }

\def\checkarrowslope #1 from #2 to #3 {{
	\ifnum\arrowslopex>4
	    \errmessage{Bad slope for arrow <#1> from #2 to #3:
			(\number\arrowslopex,\number\arrowslopey)}
	  \fi
	\ifnum\arrowslopex<-4
	    \errmessage{Bad slope for arrow <#1> from #2 to #3:
			(\number\arrowslopex,\number\arrowslopey)}
	  \fi
	\ifnum\arrowslopey>4
	    \errmessage{Bad slope for arrow <#1> from #2 to #3:
			(\number\arrowslopex,\number\arrowslopey)}
	  \fi
	\ifnum\arrowslopey<-4
	    \errmessage{Bad slope for arrow <#1> from #2 to #3:
			(\number\arrowslopex,\number\arrowslopey)}
	  \fi
    }}

\def\checksegmentslope #1 from #2 to #3 {{
	\ifnum\arrowslopex>6
	    \errmessage{Bad slope for segment <#1> from #2 to #3:
			(\number\arrowslopex,\number\arrowslopey)}
	  \fi
	\ifnum\arrowslopex<-6
	    \errmessage{Bad slope for segment <#1> from #2 to #3:
			(\number\arrowslopex,\number\arrowslopey)}
	  \fi
	\ifnum\arrowslopey>6
	    \errmessage{Bad slope for segment <#1> from #2 to #3:
			(\number\arrowslopex,\number\arrowslopey)}
	  \fi
	\ifnum\arrowslopey<-6
	    \errmessage{Bad slope for segment <#1> from #2 to #3:
			(\number\arrowslopex,\number\arrowslopey)}
	  \fi
    }}

\def\arrowlabel #1{{
	 \newpoint mid
	 \getvectormidpoint arrow into mid
	 \labelbox(\midx,\midy)(\arrowslopex,\arrowslopey){#1}
    }}

\def\arrow #1 #2 #3;{{
	\calcarrow {#1} {#2}
	\checkarrowslope {#3} from #1 to #2
	\drawvector arrow
	\arrowlabel {#3}
    }}

\def\segment #1 #2 #3;{{
	\calcarrow {#1} {#2}
	\checksegmentslope {#3} from #1 to #2
	\drawline arrow
	\arrowlabel {#3}
    }}

\def\dblarrow #1 #2 #3;{{
	\calcarrow {#1} {#2}
	\checkarrowslope {#3} from #1 to #2
	\drawdoublevector arrow
	\arrowlabel {#3}
    }}

\def\rightcell #1,#2 #3;{{
	\newpoint {}
	\copypoint coordinatescale to {}
	\multiply\x by #1
	\multiply\y by #2
	\put(\x,\y){\makebox(0,0){${\displaystyle{#3}}\atop\Longrightarrow$}}
    }}

\def\leftcell #1,#2 #3;{{
	\newpoint {}
	\copypoint coordinatescale to {}
	\multiply\x by #1
	\multiply\y by #2
	\put(\x,\y){\makebox(0,0){${\displaystyle{#3}}\atop\Longleftarrow$}}
    }}

\def\isoarrow #1 #2;{{
	\calcarrow {#1} {#2}
	\checkarrowslope {isomorphism} from #1 to #2
	\drawdoublevector arrow
	\ifnum0=\arrowslopex
	    \arrowlabel {\wr}
	  \else
	    \arrowlabel {\sim}
	  \fi
    }}

\def\calcrectlength {
	\newlocalcount\cornerlength
	\cornerlength= \cornersize
	\divide\cornerlength by \unitlength
    }

\def\calcdiaglength {
	\calcrectlength
	\multiply\cornerlength by 1000
	\divide\cornerlength by 1414
    }

\def\westcorner #1;{{
	\calcdiaglength

	\newpoint {}
	\getobjectposition {#1} into {}
	\advance\x by -\cornerlength
	\advance\x by -\cornerlength

	\addpoint nudge to {}
	\setpoint nudge to (0,0)

	\put(\x,\y){\line(1, 1){\cornerlength}}
	\put(\x,\y){\line(1,-1){\cornerlength}}
    }}

\def\eastcorner #1;{{
	\calcdiaglength

	\newpoint {}
	\getobjectposition {#1} into {}
	\advance\x by \cornerlength
	\advance\x by \cornerlength

	\addpoint nudge to {}
	\setpoint nudge to (0,0)

	\put(\x,\y){\line(-1, 1){\cornerlength}}
	\put(\x,\y){\line(-1,-1){\cornerlength}}
    }}

\def\northcorner #1;{{
	\calcdiaglength

	\newpoint {}
	\getobjectposition {#1} into {}
	\advance\y by \cornerlength
	\advance\y by \cornerlength

	\addpoint nudge to {}
	\setpoint nudge to (0,0)

	\put(\x,\y){\line( 1,-1){\cornerlength}}
	\put(\x,\y){\line(-1,-1){\cornerlength}}
    }}

\def\southcorner #1;{{
	\calcdiaglength

	\newpoint {}
	\getobjectposition {#1} into {}
	\advance\y by -\cornerlength
	\advance\y by -\cornerlength

	\addpoint nudge to {}
	\setpoint nudge to (0,0)

	\put(\x,\y){\line( 1, 1){\cornerlength}}
	\put(\x,\y){\line(-1, 1){\cornerlength}}
    }}

\def\northeastcorner #1;{{
	\calcrectlength

	\newpoint {}
	\getobjectposition {#1} into {}
	\advance\x by \cornerlength
	\advance\y by \cornerlength

	\addpoint nudge to {}
	\setpoint nudge to (0,0)

	\put(\x,\y){\line( 0,-1){\cornerlength}}
	\put(\x,\y){\line(-1, 0){\cornerlength}}
    }}

\def\northwestcorner #1;{{
	\calcrectlength

	\newpoint {}
	\getobjectposition {#1} into {}
	\advance\x by -\cornerlength
	\advance\y by \cornerlength

	\addpoint nudge to {}
	\setpoint nudge to (0,0)

	\put(\x,\y){\line( 0,-1){\cornerlength}}
	\put(\x,\y){\line( 1, 0){\cornerlength}}
    }}

\def\southeastcorner #1;{{
	\calcrectlength

	\newpoint {}
	\getobjectposition {#1} into {}
	\advance\x by \cornerlength
	\advance\y by -\cornerlength

	\addpoint nudge to {}
	\setpoint nudge to (0,0)

	\put(\x,\y){\line( 0, 1){\cornerlength}}
	\put(\x,\y){\line(-1, 0){\cornerlength}}
    }}

\def\southwestcorner #1;{{
	\calcrectlength

	\newpoint {}
	\getobjectposition {#1} into {}
	\advance\x by -\cornerlength
	\advance\y by -\cornerlength

	\addpoint nudge to {}
	\setpoint nudge to (0,0)

	\put(\x,\y){\line( 0, 1){\cornerlength}}
	\put(\x,\y){\line( 1, 0){\cornerlength}}
    }}

\long\def\fig#1{{
\newlocalcount\leastobjectx
\newlocalcount\leastobjecty
\newlocalcount\greatestobjectx
\newlocalcount\greatestobjecty
\leastobjectx= 2147483647
\leastobjecty= 2147483647
\greatestobjectx= -2147483647
\greatestobjecty= -2147483647
\newlocalbox\thefigurebox
\savebox{\thefigurebox}{#1}
\advance\greatestobjectx by -\leastobjectx
\advance\greatestobjecty by -\leastobjecty
\begin{center}
\begin{picture}(\greatestobjectx,\greatestobjecty)
	       (\leastobjectx,\leastobjecty)
\put(0,0){\usebox{\thefigurebox}}
\end{picture}
\end{center}
}}

\usepackage{amsmath}
\usepackage{graphicx,psfrag}  
\usepackage[psamsfonts]{amssymb}
\usepackage{amscd}
\usepackage{floatflt}
 
\title{On the ergodicity of partially hyperbolic systems}
\author{Keith Burns and Amie Wilkinson}
\newtheorem{theorem}{Theorem}[section]
\newtheorem{proposition}[theorem]{Proposition}
\newtheorem{lemma}[theorem]{Lemma}
\newtheorem{corollary}[theorem]{Corollary}

\def\proof{{\bf {\medskip}{\noindent}Proof. }}

\def\remark{{\bf {\bigskip}{\noindent}Remark: }}

\def\endproof{$\diamond$ \bigskip}

\def\title{\em}
\def\eps{\varepsilon}

\def\W{\mathcal{W}}
\def\hW{{\widehat{\mathcal{W}}}}
\def\hE{{\widehat{E}}}
\def\hm{{\hat{m}}}

\def\P{\mathcal{P}}
\def\F{\mathcal{F}}

\def\hB{\widehat{B}}
\def\hJ{\widehat{J}}
\def\transverse{\,\raise2pt\hbox to1em{\hfil$\top$\hfil}\hskip -1em \hbox
to1em{\hfil$\cap$\hfil}\,} 
\newcommand\R{\mbox{\bf R}}

\newlength{\figboxwidth} 

\begin{document}
\input{epsf.tex}
\maketitle
\begin{abstract} Pugh and Shub \cite{PS} have conjectured that essential accessibility implies ergodicity,
for a $C^2$, partially hyperbolic, volume-preserving diffeomorphism. We prove this conjecture under a mild center bunching assumption, which is satsified by
all partially hyperbolic systems with $1$-dimensional center bundle.  We also obtain ergodicity results for $C^{1+\gamma}$ partially hyperbolic systems.

\end{abstract} 
\section*{Introduction}

In \cite{hopf} Eberhard Hopf introduced a simple argument that
proved the ergodicity (with respect to Liouville measure)
of the geodesic flow of a compact, negatively curved
surface.  The argument has since been applied to increasingly general classes of dynamical systems.  The key feature that these systems possess is 
hyperbolicity.  The strongest form of hyperbolicity is uniform hyperbolicity.  A diffeomorphism $f:M\to M$ of a compact manifold $M$ is {\em uniformly hyperbolic} or {\em Anosov} if there exists a splitting of the tangent 
bundle into $Tf$-invariant subbundles:
$$TM = E^s\oplus E^u,$$
and a continuous Riemannian metric, such that for every unit vector
$v\in TM$:
\begin{alignat}{4}\label{e=defn0}
\|&Tf v\| < 1 &\qquad \hbox{if }v\in E^s, \\
\|&Tf v\|  > 1 &\qquad \hbox{if } v\in E^u.
\end{alignat}
Anosov flows are defined similarly, with $E^s\oplus E^u$ complementary
to the bundle $E^0$ that is tangent to the flow direction. 
The bundles $E^s$ and $E^u$ of an Anosov system
are tangent to the {\em stable} and {\em unstable} 
foliations $\W^s$ and $\W^u$,
respectively.  The properties of these foliations play a crucial role in
the Hopf argument.

Hopf's original argument established ergodicity for volume-preserving
uniformly hyperbolic systems under the assumption that the foliations  $\W^s$ and $\W^u$
are $C^1$.  While the leaves of these foliations are always as smooth as
the diffeomorphism, the foliations are usually only only H\"older 
continuous in the direction transverse to the leaves.    
In particular, for geodesic flows on arbitrary compact manifolds
of negative sectional curvature, these foliations are not always $C^1$.
 
Anosov and Sinai \cite{AS, An} observed that  the $C^1$ condition on the
stable and unstable foliations in the Hopf argument
could be replaced by the weaker condition of absolute continuity,
which we discuss in Section~\ref{s=abscont}.
They showed that $\W^s$ and $\W^u$ are absolutely continuous
if the system is  $C^2$, thereby establishing
ergodicity of all $C^2$ volume-preserving uniformly hyperbolic systems, 
including geodesic flows for compact manifolds 
of negative sectional curvature.  

At this point, it became clear that the Hopf argument should extend to 
even more general settings.  Two natural generalizations of 
uniform 
hyperbolicity are:
\begin{itemize}
 \item nonuniform hyperbolicity, which requires hyperbolicity along almost every orbit, but allows the expansion of $E^u$ and the contraction of $E^s$ to
weaken near the exceptional set where there is no hyperbolicity; and
 \item partial hyperbolicity, which requires uniform expansion of $E^u$ and uniform contraction of $E^s$, but allows central directions at each point, in which the expansion and contraction is dominated by the behavior in the hyperbolic directions.
\end{itemize}
The first direction is Pesin theory; the second is the subject of this paper.

Brin and Pesin \cite{BP1} and independently Pugh and Shub \cite{PSanosovactions}
first examined the ergodic properties of partially hyperbolic systems soon
after the work of Anosov and Sinai. The current definition
of partial hyperbolicity is more general than theirs, but has the
same basic features.\footnote{The difference is that in \cite{BP1} the functions
$\nu,\hat\nu,\gamma$ and $\hat\gamma$ in the definition of partial hyperbolicity are assumed to be constant.} 
We say that a diffeomorphism 
$f:M\to M$ of a compact manifold $M$ is partially hyperbolic
if the following conditions hold.
There is a nontrivial splitting of the tangent bundle,
$TM=E^s\oplus E^c\oplus E^u$, that is invariant under
the derivative map $Tf$.
Further, there is
a Riemannian metric for which we can choose continuous positive 
functions $\nu$,  $\hat\nu$, $\gamma$ and $\hat\gamma$ with  
\begin{eqnarray}\label{e=defn}
 \nu, \hat\nu < 1\quad\hbox{ and } \quad
 \nu <  \gamma < \hat\gamma^{-1} < \hat\nu^{-1} 
 \end{eqnarray}
such that, for any unit vector $v\in T_pM$,
\begin{alignat}{4}\label{e=defn2}
\|&Tf v\| <  \nu(p) ,&\qquad \hbox{if }v\in E^s(p),  \\
\gamma(p) < \|&Tf v\|  < 
\hat\gamma(p)^{-1} ,&\qquad \hbox{if } v\in E^c(p), \\
\hat\nu(p)^{-1}< \|&Tf v\|  ,&\qquad \hbox{if } v\in E^u(p).
\end{alignat}
Partial hyperbolicity is a $C^1$-open condition: any diffeomorphism sufficiently $C^1$-close to a partially hyperbolic diffeomorphism is itself partially hyperbolic. Partially hyperbolic flows are defined similarly. For an extensive discussion of examples of partially hyperbolic dynamical systems, see the survey article \cite{BPSW} and the book \cite{P}. Among these examples are: the time-$1$ map of an Anosov flow, the frame flow for a compact manifold of negative sectional curvature, and many affine transformations of compact homogeneous spaces.

As in the Anosov case, the stable and unstable bundles $E^s$ and $E^u$ of
a partially hyperbolic diffeomorphism are tangent to foliations, 
which we again denote by $\W^s$ and $\W^u$ respectively \cite{BP1}. 
Brin-Pesin and Pugh-Shub proved that these foliations 
are absolutely continuous.  

By its very nature, the Hopf argument shows that 
for almost every $p\in M$, almost every point of $\W^s(p)$ and 
almost every point of $\W^u(p)$ lies in the ergodic component of $p$ (cf. \cite{xia}).  Thus we can only hope to prove ergodicity using a Hopf argument
if something close to the following condition holds.

\begin{dfinition}
A partially hyperbolic diffeomorphism $f:M\to M$ is {\em accessible} if
any point in $M$ can be reached from any other along an {\em $su$-path}, which
is a concatenation of finitely many subpaths, 
each of which lies entirely in a single leaf of $\W^s$ or a single leaf of $\W^u$.
\end{dfinition}

\medskip

The {\em accessibility class of $p\in M$} is the set of all $q\in M$ that can be reached from $p$ along an $su$-path.  Accessibility means that there is one accessibility class, which contains all points.  The following notion is a natural weakening of
accessibility.

\begin{dfinition}
A partially hyperbolic diffeomorphism $f:M\to M$ is {\em essentially
accessible} if every measurable set that is a union of entire accessibility classes has either full or zero measure.
\end{dfinition}

\medskip

Pugh and Shub have conjectured that essential accessibility implies ergodicity,
for a $C^2$, partially hyperbolic, volume-preserving diffeomorphism \cite{PSa}. We prove this conjecture under one, rather mild additional assumption.

\begin{dfinition}
A partially hyperbolic diffeomorphism is {\em center bunched} if 
 the functions
$\nu, \hat\nu, \gamma$, and $\hat\gamma$
can be chosen so that:
\begin{eqnarray}\label{e=bunch}
  \nu < \gamma\hat\gamma  \qquad\hbox{and}\qquad    \hat\nu < \gamma\hat\gamma.
\end{eqnarray}
\end{dfinition}

Our main result is:

\begin{theorem}\label{t=main} Let $f$ be $C^2$, volume preserving,
partially hyperbolic,  and center bunched.  If $f$ is essentially accessible, then $f$
is ergodic, and in fact has the Kolmogorov property.
\end{theorem}

This result extends earlier results about ergodicity of partially hyperbolic 
systems. Brin and Pesin \cite{BP1}  
proved in the early 1970's
that a $C^2$ volume-preserving partially hyperbolic diffeomorphism 
that is essentially accessible 
is ergodic if it satisfies the following additional conditions:
\begin{itemize}
\item {\em Center bunching:} Inequalities (\ref{e=bunch}) hold. 
\item {\em Dynamical coherence:} There are foliations $\W^c$, $\W^{cs}$ and
$\W^{cu}$ tangent to $E^c$, $E^c \oplus E^s$ and $E^c \oplus E^u$ respectively.
\item {\em Lipschitzness of $\W^c$:} There are Lipschitz foliation charts for 
$\W^c$.
\end{itemize}
While the Brin-Pesin argument applies to many examples of partially hyperbolic
diffeomorphisms, the third condition is in some ways very restrictive:  
Lipschitzness of 
$\W^c$ can be destroyed by arbitrarily small perturbations \cite{SW2}.
Brin and Pesin's theorem applies in particular 
to  the time-$1$ map $\varphi_1$ of the geodesic flow for 
a compact surface of constant negative curvature.  
If we make a $C^1$ small perturbation
to $\varphi_1$, all of Brin and Pesin's hypotheses continue to hold, except 
Lipschitzness of $\W^c$ \cite{D}.

It was not until the 1990's that Grayson, Pugh and Shub \cite{GPS}
were able to show that
any small perturbation of $\varphi_1$ is ergodic; in other words, $\varphi_1$
is {\em stably ergodic:} any $C^2$ volume-preserving diffeomorphism 
sufficiently $C^1$-close to $\varphi_1$ is ergodic. The ideas in this groundbreaking paper
have been generalized in several stages \cite{amiethesis, PSa, PS}, culminating
in  \cite{PS}.  The main result of \cite{PS} assumes dynamical coherence and
uses a significantly stronger version of center bunching than inequalities
(\ref{e=bunch}).  The center bunching hypothesis 
in \cite{PS} requires that the action of $Tf$ on $E^c$ be close to 
{\em isometric} --- that is, both $\gamma$ and $\hat\gamma$ (and not just their
product) must be close to $1$.

%{\bf Somewhere, maybe not here discuss the difference between our center bunching assumption and theirs.}

By contrast, our center bunching hypothesis requires only that
 the action of $Tf$ on $E^c$ be close enough to 
{\em conformal} that the hyperbolicity of $f$ 
dominates the nonconformality of $Tf$ on $E^c$. 
Center bunching always holds
when $Tf\vert_{E^c}$ is conformal. For then we have $\|T_pf v\| = \|T_pf|_{E^c(p)}\|$
for any unit vector $v \in E^c(p)$, and hence we can choose $\gamma(p)$ 
slightly  smaller and $\hat\gamma(p)^{-1}$ slightly bigger than
 $$
 \|T_pf|_{E^c(p)}\|.
 $$
By doing this we may make the ratio $\gamma(p)/\hat\gamma(p)^{-1} = \gamma(p)\hat\gamma(p)$ arbitrarily close to 1, and hence larger than both $\nu(p)$ and
$\hat\nu(p)$.

In particular, center bunching holds whenever $E^c$ is one-dimensional.
As a corollary, we obtain:

\begin{corollary}\label{c=dim1}
Let $f$ be $C^{2}$, volume preserving and
partially hyperbolic with $\hbox{\rm dim}(E^c) = 1$.  If $f$ is essentially accessible, then $f$
is ergodic, and in fact has the Kolmogorov property. 
\end{corollary}

This establishes the Pugh-Shub Conjecture mentioned above in the case
where $E^c$ is 1-dimensional.

Corollary~\ref{c=dim1} has also been recently proved 
by  F.  Rodr\'iguez Hertz,  J. Rodr\'iguez Hertz, and R. Ures \cite{HHU}. 
Their  argument is mainly based on techniques 
in an earlier version\footnote{This earlier version proved the same result as the present paper but under the additional hypothesis of dynamical coherence, i.e. the existence of foliations tangent to the bundles $E^c \oplus E^u$ and  $E^c \oplus E^s$.} 
of the present paper \cite{BW}.
They prove in addition that stable accesibility is $C^r$ dense
among the $C^r$ partially hyperbolic diffeomorphisms with $1$-dimensional
center, which implies, for $r\geq 2$, that stable ergodicity is $C^r$ dense
among the volume-preserving 
$C^r$ partially hyperbolic diffeomorphisms with one-dimensional 
center.  Their work establishes
the main stable ergodicity conjectures of Pugh and Shub 
(\cite{PS}, Conjectures 1-3) in the case where $E^c$ is one-dimensional.

There is only one place in the proof of Theorem~\ref{t=main} where we need the diffeomorphism
to be $C^2$ as opposed to $C^{1+\delta}$.  This is when we use the fact that center bunching
implies that the stable and unstable holonomies between center leaves are
Lipschitz.   This fact is proved using a graph transform argument in
\cite{BP1} and also (in a slightly more general setting) in \cite{PSW, PSWc}.
In \cite{BWhold}, we show that the same result about holonomies 
holds when $C^2$ is replaced by $C^{1+\delta}$, at the expense of a 
more stringent bunching hypothesis.  Plugging this result into
the proof of Theorem~\ref{t=main}, we obtain:

\begin{theorem}\label{t=mainholder} Let $f$ be $C^{1+\delta}$, volume preserving,
and partially hyperbolic. Let $\mu,\hat \mu$ be continuous functions
satisfying:
\begin{alignat}{4}\label{e=defn3}
\mu(p) < \|&Tf v\| ,&\qquad \hbox{if }v\in E^s(p),  \\
\|&Tf v\| < \hat\mu(p)^{-1} ,&\qquad \hbox{if } v\in E^u(p).
\end{alignat}
Suppose that $f$ satisfies the {\em strong center bunching condition}:
\begin{eqnarray}\label{e=strongbunch}
\nu^{\theta}  < \gamma\hat\gamma \qquad\hbox{and}\qquad {\hat\nu}^{\theta}  < \gamma\hat\gamma,
\end{eqnarray}
where ${\theta}\in(0,\delta)$  satisfies the inequalities:
\begin{eqnarray}\label{e=exp}
\nu\gamma^{-1} < \mu^{\theta},\qquad \hat\nu\hat\gamma^{-1} < \hat\mu^{\theta}.
\end{eqnarray}

If $f$ is essentially accessible, then $f$
is ergodic, and in fact has the Kolmogorov property.
\end{theorem}

We remark that any $\theta$ satisfying the conditions in (\ref{e=exp}) is
a H\"older exponent for the central distribution $E^c$ (see, e.g.,
Theorem A in \cite{PSW}). It would be interesting to know whether 
(\ref{e=strongbunch}) could be replaced by (\ref{e=bunch})
in Theorem~\ref{t=mainholder}. The strong center bunching condition in (\ref{e=strongbunch}) is automatically satisfied when
$E^c$ is one-dimensional, since, as above,
we may then choose $\gamma$ and $\hat\gamma^{-1}$
so that $\gamma\hat\gamma$ is arbitrarily close to $1$.
As a corollary of Theorem~\ref{t=mainholder}, we obtain:

\begin{corollary}\label{c=dim1holder}
Let $f$ be $C^{1+\delta}$, volume preserving and
partially hyperbolic with $\hbox{\rm dim}(E^c) = 1$.  If $f$ is essentially accessible, then $f$
is ergodic, and in fact has the Kolmogorov property. 
\end{corollary}

We do not know whether the bunching assumption can 
be dropped in 
Theorem~\ref{t=main}.   A first step in answering this question might be to 
answer the following question:

\begin{question} Suppose that $f$ is $C^2$, volume preserving 
and partially hyperbolic.  If $E^c$ is absolutely continuous 
(Lipschitz, smooth...) and $f$ is (essentially)
accessibile, is $f$ then ergodic?
\end{question}

\bigskip

As a concrete example, consider a diffeomorphism 
$f_\lambda:{\bf T}^2\times{\bf T}^2 \to
{\bf T}^2\times{\bf T}^2$ of the form:
$$f_\lambda (x,y) = (A(x), g_\lambda(y)),$$
where $A:{\bf T}^2\to {\bf T}^2$ is the linear Anosov diffeomorphism
given by 
$$A = \left(\begin{array}{cc}
2 & 1 \\
1 & 1
\end{array}
\right)^2,$$
and $g_\lambda: {\bf T}^2\to {\bf T}^2$ is a standard map of the form:
$$g_\lambda(z,w) = (z + w, w + \frac{\lambda}{2\pi}\sin(2\pi(z+w))).$$
A straightforward calculation shows that
there is an interval $\Lambda\subset {\bf R}$ containing $(-4,4)$  
such that, if $\lambda\in \Lambda$, then  $f_\lambda$ is partially
hyperbolic with respect to the standard (flat) metric on
${\bf T}^2\times {\bf T}^2$.   It is also not difficult to show,
by examining the spectrum of $Tf_\lambda$ at the fixed point $(0,0)$, that
$f_\lambda$ is center bunched if and only if $\lambda\in (-1,1)$.

This example appears in \cite{SW1}, where it is shown that 
there is a function $\varphi:{\bf T}^2\to {\bf T}^2$ with
$\varphi(0) = 0$
and an interval $E = (0,\epsilon_0)$
such that, for all $(\lambda,\epsilon)\in \Lambda\times E$, the map
$$f_{\lambda,\epsilon}(x,y) = (A(x), g_\lambda(y) +\epsilon\varphi(x))$$
is both partially hyperbolic and stably accessible. Futhermore,
$f_{\lambda,\epsilon}$ is center bunched if and only if 
$(\lambda,\epsilon)\in (-1,1)\times E$. For all $(\lambda,\epsilon)\in \Lambda\times E$, the center bundle
$E^c$ is tangent to the fibers $\{x\}\times {\bf T}^2$ and is $C^\infty$.
The foliations $\W^{cu}$ and $\W^{cs}$ are also $C^\infty$.

Theorem~\ref{t=main} implies that $f_{\lambda,\epsilon}$ is
stably ergodic for all $(\lambda,\epsilon)\in (-1,1)\times E$.
We do not however know whether $f_{\lambda,\epsilon}$ is
ergodic for even a single value of $(\lambda,\epsilon)\in \left(\Lambda\setminus (-1,1)\right)\times E$.

\bigskip

Theorem~\ref{t=main} is proved in Section~\ref{s=mainthm} as a consequence of 
Theorem~\ref{t=reform}, which is really the central result of the paper. 
Theorem~\ref{t=reform} is proved in Sections~\ref{s=jqconf} and \ref{s=density}.
We thank Marcelo Viana for telling
us about the proof of absolute continuity of stable foliations in
the pointwise partially hyperbolic setting in \cite{viana}.
We thank Charles Pugh and Mike Shub for very useful comments.  
Keith Burns was supported by
NSF grants DMS-0100416 and DMS-0408704, and Amie Wilkinson 
by NSF grants DMS-0100314
and DMS-0401326.
A proof of Theorem~\ref{t=main} under the additional hypothesis 
of dynamical coherence appeared earlier in the unpublished preprint \cite{BW}.

%\section{Outline of the proof of the main result}

\section{Preliminaries}

\subsection{Notational conventions}

We use the  convention that if
 $q$ is a point in $M$ and $j$
is an integer, then $q_j$ denotes the point $f^j(q)$, with $q_0=q$.
If $\alpha:M\to {\bf R}$ is a positive function, and $j\geq 1$ is
an integer, let
$$
\alpha_j(p) = \alpha(p)\alpha(p_1)\cdots\alpha(p_{j-1}), 
$$
and
$$
\alpha_{-j}(p) = \alpha(p_{-j})^{-1}\alpha(p_{-j+1})^{-1}\cdots\alpha(p_{-1})^{-1}.
$$
 We set $\alpha_0(p) = 1$.
Observe that $\alpha_j$ is a
multiplicative cocycle; in particular, we have
$\alpha_{-j}(p)^{-1}= \alpha_j(p_{-j})$. Note also that $(\alpha\beta)_j = \alpha_j\beta_j$, and if 
$\alpha$ is a constant function, then $\alpha_n = \alpha^n$.

The notation $\alpha < \beta$, where 
$\alpha$ and $\beta$ are continuous functions, means that the inequality
holds pointwise, and the function $\min\{\alpha,\beta\}$ takes
the value $\min\{\alpha(p),\beta(p)\}$ at the point $p$.

As usual $P = O(Q)$ means that there is
a constant $C > 0$ such that $|P| \leq CQ$. 
Usually $P$ and $Q$ will depend on an integer
$n$ and one or more points in $M$. 
The constant $C$ must be independent of $n$ and the choice of the points.

\subsection{Foliation boxes and local leaves}\label{s=local}

Let $\F$ be a foliation of an $n$-manifold $M$ 
with $d$-dimensional smooth leaves. 
For $r>0$, we denote by 
$\F(x,r)$ the connected component of $x$ 
in the intersection of $\F(x)$ with the ball $B(x,r)$.
%For $B$ any subset of $M$, we set:
%$$\F(B,r) = \bigcup_{x\in B} \F(x,r).$$

A {\em foliation box for
$\F$} is the image $U$ of 
$\R^{n-d}\times \R^{d}$ under a homeomorphism that sends each
vertical $\R^d$-slice into a leaf of $\F$. The images of the
vertical $\R^d$-slices will be called {\em local leaves of $\F$ in $U$}.

A {\em smooth transversal} to $\F$ in $U$ is a smooth codimension-$d$
disk in $U$ that intersects each local leaf in $U$ exactly once and whose
tangent bundle is uniformly transverse to $T\F$.
If $\tau_1$ and $\tau_2$ are two smooth transversals to $\F$ in $U$,
we have the {\em holonomy map} $h_{\F}: \tau_1 \to \tau_2$,
which takes a point in $\tau_1$ to the intersection  of its local leaf
in $U$ with $\tau_2$.

\subsection{Adapted metrics}\label{s=adapted}
We assume that the Riemannian metric on $M$ is chosen so that the inequalities
  (\ref{e=defn})--(6) % \ref{e=bunch}
involving $\nu,\gamma, \hat\nu,\hat\gamma$ in the Introduction hold.  
Such a metric will be called 
adapted.  Note that a rescaling of an adapted metric is still adapted. It will be
convenient to assume that the metric is scaled so that the geodesic balls of radius $1$ are
very small neighborhoods of their centers. Distance with respect to the metric will
be denoted by $d$.  

There is no harm in increasing $\nu$ and
$\hat\nu$ and decreasing $\gamma, \hat\gamma$  slightly,
provided that the inequalities 
     (\ref{e=defn})--(6)  %\ref{e=bunch})
 still hold.  If $f$ is center bunched, 
the change must also be small enough so that
inequalities (\ref{e=bunch}) still hold.
Similarly, if  $f$ is strongly center bunched, 
the change must also be small enough so that
inequalities (\ref{e=strongbunch}) still hold.
%We will assume that the functions  
%$\nu,\hat\nu,\gamma,\hat\gamma$ are smooth, 
%though 
%H\"older continuity will suffice for our purposes.

By rescaling the metric
on $M$, we may assume that for some $R > 1$, and any $x \in M$, 
the Riemannian ball $B(x,R)$ is contained in
foliation boxes for both $\W^s$ and $\W^u$.
We assume that $R$ is large enough so that 
all the objects considered in the sequel
are small compared with $R$.
Having fixed such an $R$, we define, for
$a= s$ or $u$, the {\em local leaf of $\W^a$ through $x$}
by: 
$$
\W^a_{loc}(x) = \W^a(x,R).
$$ 
Any foliation box $U$ for either  $\W^s$ or $\W^u$
that we consider in the rest of the paper will be small enough so that 
$\W^a_{loc}(x)\cap U$ is a local leaf of $\W^a$ in $U$ for each $x\in U$.
By (if necessary)
further rescaling the metric to make the local leaves smaller, 
we may assume that our metric is still adapted, and that
for all $p\in M$, and $q,q'\in B(p,R)$, 
we have the following:
\begin{eqnarray}\label{e=locdist}
q\in \W^s_{loc}(q')\quad \Longrightarrow\quad
d(f(q),f(q')) \leq \nu(p) d(q,q'),
\end{eqnarray}
and similarly,
\begin{eqnarray*}
q\in \W^u_{loc}(q')\quad \Longrightarrow\quad
d(f^{-1}(q),f^{-1}(q')) \leq \hat\nu(f^{-1}(p)) d(q,q').
\end{eqnarray*}
In particular, $f(\W^s_{loc}(p)) \subset \W^s_{loc}(f(p))$ and
$f^{-1}(\W^u_{loc}(p)) \subset \W^u_{loc}(f^{-1}(p))$, for all $p\in M$.
This is possible because $\nu$ and $\hat\nu$ are continuous, and
the inequalities (4) and (6) that they satisfy are strict.

An inductive argument then gives:
\begin{lemma}\label{l=distest} 
If $q_j,q_j'\in B(p_j,R)$ for
$j=0,\ldots, n-1$, and
$q \in \W^s_{loc}(q')$, then
$$
d(q_n,q_n') \leq \nu_n(p) d(q,q').
$$

If $q_{-j},q_{-j}'\in B(p_{-j},R)$ for
$j=0,\ldots, n-1$, and
$q \in \W^u_{loc}(q')$, then
$$d(q_{-n},q_{-n}') \leq \hat\nu_{-n}(p)^{-1} d(q,q').$$
\end{lemma}
\proof We prove the first claim; the second follows from 
the first, with $f$ replaced by $f^{-1}$. The proof is by
induction on $n$.  The claim is vacuously true for $n=0$.  Suppose
the claim holds for $n=k$.  The inductive assumption gives that
$$
d(q_k,q_k') \leq \nu_k(p) d(q,q')
$$
and $q_k'\in \W^s_{loc}(q_k)$.  Then (\ref{e=locdist}),
applied at $p_k$, implies that  $q_{k+1}'\in \W^s_{loc}(q_{k+1})$,
and
\begin{eqnarray*}
d(q_{k+1},q_{k+1}') &\leq& \nu(p_k) d(q_k,q_k')\\
&\leq& \nu(p_k) \nu_k(p) d(q,q')\\
&=& \nu(p_{k+1}) d(q,q').
\end{eqnarray*}
\endproof

\subsection{Volume and density}

When we say that the diffeomorphism $f$ is
volume preserving, we mean that
$f$ preserves a measure $m$ that
is equivalent to the Riemannian volume $m_M$ on $M$. (This definition is
independent of the metric, since the volumes defined by two different metrics 
on $M$ are always equivalent.)
Unless otherwise specified,
measurable will mean measurable with respect to $m$.

If $S \subseteq M$ is a smooth submanifold, we denote by 
$m_S$ the volume of the 
induced Riemannian metric on $S$. If $\F$ is a foliation with smooth leaves,
and $A$ is contained in a single leaf of $\F$ and is measurable in that 
leaf, then we denote by $m_\F(A)$ the induced Riemannian 
volume of $A$ in that leaf. 
A  set is said to be 
{\em saturated by a foliation $\F$} or {\em $\F$-saturated} if it
is a union of entire leaves of $\F$.  A set $A$ is 
{\em essentially $\F$-saturated} if there exists a measurable $\F$-saturated
set $A'$, which we call an {\em essential
$\F$-saturate} of $A$,  with $m(A\,\Delta\, A') = 0$.

If $\mu$ is a measure and 
$A$ and $B$ are $\mu$-measurable sets with $\mu(B) >0$,
we define the {\em density of $A$ in $B$} by: 
 $$
 \mu(A:B) = \frac{\mu(A\cap B)}{\mu(B)}.
 $$ 
A point $x\in M$ is a {\em Lebesgue density point} 
of a measurable set $X\subseteq M$ if
$$\lim_{r\to 0} m(X: B_r(x)) = 1.$$
The Lebesgue Density Theorem implies that if $A$ is a measurable
set and $\widehat A$ is the set of Lebesgue density points of $A$,
then $m(A\,\Delta\, \widehat A) = 0$.

Lebesgue density points can be characterized using nested sequences
of measureable sets.
We say that a sequence of measurable sets $Y_n$ {\em nests} at point $x$
if $Y_0\supset Y_1 \supset Y_2 \supset \cdots \supset \{x\}$, and
$$\bigcap_n Y_n = \{x\}.$$ 
A sequence $Y_n$ that nests at $x$ is a 
{\em Lebesgue density sequence at $x$} if, for every measurable set $X$, 
$x$ is a Lebesgue density point of $X$ if and only if:
$$\lim_{n\to\infty} m(X:Y_n) = 1.$$
It is easily shown that a Lebesgue density sequence $Y_n$ 
must be regular.

\begin{dfinition} A sequence of measurable sets $Y_n$ is {\em regular} if
 there exist $\delta> 0$ and $l\geq 1$ such that, for all $n\geq 0$, 
$$m(Y_{n+l})\geq \delta m(Y_n).$$
\end{dfinition}
The simplest example of a Lebesgue density sequence at $x$ is the
sequence of balls $B(x,\rho^n)$, where $\rho\in (0,1)$.

In our proof of Theorem~\ref{t=main}, we characterize the Lebesgue density
points of a special class of measurable sets, those that are both 
{\em essentially
$\W^s$-saturated and essentially $\W^u$-saturated}.  
We say that $Y_n$ is an {\em $su$-density sequence at $x$} if
$Y_n$ nests at $x$, $Y_n$ is regular,
and, for every set $X$ that is both essentially
$\W^s$-saturated and essentially $\W^u$-saturated, 
$x$ is a Lebesgue density point of $X$ if and only if:
$$\lim_{n\to\infty} m(X:Y_n) = 1.$$
Note that an $su$-density sequence is not necessarily a Lebesgue density
sequence, because only certain measurable sets are considered in the definition of 
an $su$-density sequence. In fact, many of the $su$-density sequences constructed in
this paper are not Lebesgue density sequences.

In our proof we will frequently have to pass the property of being 
an $su$-density sequence from
one sequence $Y_n$ that nests at $x$ to another sequence $Z_n$
that nests at $x$. In order to do so, we have to 
show that $Z_n$ is also regular and that, for every measurable set $X$
that is both essentially
$\W^s$-saturated and essentially-$\W^u$ saturated,
\begin{eqnarray}\label{e=denseq}
\lim_{n\to\infty} m(X:Y_n) = 1 \quad\Longleftrightarrow \quad \lim_{n\to\infty} m(X:Z_n) = 1.
\end{eqnarray}
We have two techniques for doing this.
Which technique we use depends on the construction of $Y_n$ and $Z_n$.  

The first technique is very simple. Two sequences of sets $Y_n$ and $Z_n$ are {\em comparable} if there exists
a $k\geq 1$ such that, for all $n\geq 0$, we have
$$Y_{n+k} \subseteq Z_n,\qquad \hbox{ and } \qquad Z_{n+k}\subseteq Y_{n}.$$
Comparability is an equivalence relation.
The following lemma is a straightforward consequence of the definitions.

\begin{lemma}\label{l=compreg}
Let $Y_n$ and $Z_n$ be comparable sequences of measurable sets, with
$Y_n$ regular.
Then $Z_n$ is also regular.  If the sets $Y_n$ have positive measure, then
so do the $Z_n$, and, for any measurable set $X$,
$$\lim_{n\to\infty} m(X:Y_n) = 1 \quad\Longleftrightarrow \quad \lim_{n\to\infty} m(X:Z_n) = 1.$$ 
\end{lemma} 

\begin{corollary}\label{c=compreg}
Suppose that $Y_n$ and $Z_n$ both nest at $x$ and are comparable.
Then $Y_n$ is an $su$-density sequence if and only if $Z_n$ is
an $su$-density sequence. 
\end{corollary}

The second technique uses absolute continuity of the foliations $\W^s$ and $\W^u$, plus the saturation properties of $X$, and is developed in the next
subsection.

\subsection{Absolute continuity and saturated sets}\label{s=abscont}

Our arguments in this paper use two versions of the property 
of absolute continuity of a foliation.

The first version of absolute continuity involves
holonomy maps between transversals.  
A foliation $\F$ with smooth leaves is {\em transversely 
absolutely continuous with bounded 
Jacobians} if  for every angle $\alpha\in(0,\pi/2]$,
there exists $C\geq 1$ such that, for every foliation
box $U$ of diameter less than $R$, any two smooth 
transversals $\tau_1, \tau_2$ to
$\F$ in $U$ of angle at least $\alpha$ with $\F$, and any
$m_{\tau_1}$--measurable set $A$ contained in $\tau_1$:
\begin{eqnarray}\label{e=acholon}
 C^{-1} m_{\tau_1}(A) \leq m_{\tau_2}(h_\F(A))\leq C m_{\tau_1}(A).
\end{eqnarray}

The second version involves a Fubini-like property.  
A foliation $\F$ with smooth leaves is {\em absolutely continuous with 
bounded Jacobians} if, for every $\alpha\in(0,\pi/2]$, 
there exists $C\geq 1$ such that, for every 
foliation box $U$ of diameter less than $R$, any smooth 
transversal $\tau$ to
$\F$ in $U$ of angle at least $\alpha$ with $\F$, and any
measurable set $A$ contained in $U$, we have the inequality:
		
\begin{eqnarray}\label{e=fubini}
C^{-1} m(A) \leq \int_{\tau} m_{\F}(A\cap \F_{loc}(x))\, dm_\tau(x)\leq C m(A).
\end{eqnarray}

If $\F$ is transversely absolutely continuous with bounded Jacobians, then 
it is absolutely continuous with bounded Jacobians (see \cite{brinstuck} for a 
proof), 
but the converse does not hold (see Remark 3.9 in \cite{Mishapp}).  
Note that the minimal $C$ for which 
(\ref{e=acholon}) holds is not necessarily the same minimal  $C$ for which 
(\ref{e=fubini}) holds. 

The foliations $\W^s$ and $\W^u$ for a partially hyperbolic diffeomorphism
are transversely 
absolutely continuous with bounded 
Jacobians.  This was shown in the Anosov case by Anosov \cite{An},
and in the case of  partial
hyperbolicity by Brin-Pesin and Pugh-Shub \cite{BP1, PSanosovactions}.
Their proofs were written under the assumption that the function $\nu$, $\hat\nu$, 
$\gamma$ and $\hat\gamma$ are constant.
In the general case of partial hyperbolicity, where these functions are not constant,
absolute continuity
of $\W^s$ and $\W^u$ 
follows from Pesin theory.  A direct proof in this context has been
given by Abdenur and Viana \cite{viana}. 
All of these results show that the Jacobians are continuous functions,
and so are bounded, since $M$ is compact.
 In general, $\W^c$ does not have either absolute continuity
property, even when
$f$ is dynamically coherent (examples were first constructed by Katok 
\cite{Mi};
open sets of examples by Shub-Wilkinson \cite{SW2}).

The second technique mentioned in the previous subsection involves 
decomposing the sets in
a sequence nesting at $x$ along leaves of an absolutely continuous foliation.

Let $\F$ be an absolutely continuous foliation 
and let $U$ be a foliation box for $\F$. 
Let  $\tau$ be a smooth transversal to $\F$ in $U$.
Let $Y\subseteq U$ be a measurable set.  
For a point $q \in \tau$, we define the 
{\em fiber $Y(q)$ of $Y$ over $q$} to be the
intersection of $Y$ with the local leaf of $\F$ in $U$ containing $q$. 
The {\em base $\tau_Y$ of $Y$} is
the set of all $q\in \tau$ such that the fiber $Y(q)$ is
$m_\F$-measurable and  $m_\F(Y(q))>0$.
The absolute continuity of $\F$ implies that $\tau_Y$ is 
 $m_{\tau}$-measurable.
We say  ``$Y$ fibers over $Z$'' to indicate that $Z=\tau_Y$.

If, for some $c\geq 1$, the inequalities
$$
c^{-1} \leq \frac{m_\F(Y(q))}{m_\F(Y(q'))} \leq c
$$
hold for all $q,q'\in \tau_Y$, then we say that {\em $Y$ has $c$-uniform fibers}. A sequence of measurable sets $Y_n$ contained in $U$ has 
{\em $c$-uniform fibers} if each set in the sequence has $c$-uniform fibers,
with $c$ independent of $n$.

Regularity of a sequence with $c$-uniform fibers can be obtained from
regularity of its fibers and bases.

\begin{proposition}\label{p=unifreg} Suppose that the foliation $\F$ is absolutely continuous
with bounded Jacobians. Let $U$ be a foliation box for $\F$,
and let  $\tau$ be a smooth transversal to $\F$ in $U$.
Let $Y_n$ be a sequence of subsets of $U$ 
with $c$-uniform fibers.  Suppose that:
\begin{enumerate}
\item there exist $\delta_1 >0$ and $k_1 \geq 1$ such that for all $n \geq 0$
$$
m_{\tau}(\tau_{Y_{n+k_1}}) \geq \delta_1 m_{\tau}(\tau_{Y_n});
$$

\item there exist $\delta_2 > 0$ and $k_2 \geq 1$ such that for all $n \geq 0$
there are points $z\in\tau_{Y_{n+k_2}}, z'\in \tau_{Y_n}$ with
$$
 m_{\F}(Y_{n+k_2}(z)) \geq \delta_2 m_{\F}(Y_n(z')).
$$
\end{enumerate}
Then $Y_n$ is regular.
\end{proposition}

\proof Let $k= k_1k_2$ and $\delta = \min\{\delta_1^{k_2},\delta_2^{k_1}\}$.
Then both 1. and 2. hold with $k_1$ and $k_2$ replaced by $k$ and
$\delta_1$ and $\delta_2$ replaced by $\delta$.
Absolute continuity of $\F$ with bounded Jacobians implies that there
exists a $C\geq 1$ such that, for all $n\geq 0$.
$$
C^{-1} m(Y_n) \leq \int_{\tau} m_{\F}({Y_n}(q))\, dm_\tau (q)\leq C m(Y_n).
$$
This, combined with uniformity of fibers, implies that
$$
m(Y_{n+k}) \geq (c C)^{-1} m_{\F}({Y_{n+k}}(z)) m_{\tau} (\tau_{Y_n}),
$$
and
$$
m(Y_{n}) \leq c C m_{\F}({Y_n}(z')) m_{\tau} (\tau_{Y_n}).
$$
Using 1. and 2., we then obtain:
$$m(Y_{n+k}) \geq (cC)^{-2} \delta^2 m(Y_n).$$
Hence, $Y_n$ is regular.\endproof

We now turn to the second technique 
mentioned above for proving an equivalence of the form (\ref{e=denseq}). 
The main result we prove in this section is:
\begin{proposition}\label{p=compmeas}
Let $\F$ be absolutely continuous with bounded Jacobians, and let $U$ be
a foliation box for $\F$ with smooth transversal $\tau$.  
Let $\{Y_n\}$ and $\{Z_n\}$ be sequences of 
measurable
subsets of $U$ with $c$-uniform fibers, for some $c\geq 1$.
Suppose that $\tau_{Y_n} = \tau_{Z_n}$, for all $n$.  
Then, for any essentially $\F$-saturated set $X\subseteq U$, we have the equivalence:
$$
\lim_{n\to\infty} m(X:Y_
n) = 1\,\Longleftrightarrow\, 
\lim_{n\to\infty}m(X:Z_n) = 1.
$$
\end{proposition}

\begin{corollary}
Let $Y_n$ and $Z_n$ be as in Proposition~\ref{p=compmeas},
with $\F = \W^s$ or $\W^u$,
and suppose that $Y_n$ and $Z_n$ both nest at $x$.
If $Y_n$ is an $su$-density sequence at $x$, and $Z_n$ is regular, then
$Z_n$ is an $su$-density sequence at $x$.
\end{corollary}

Before proving Proposition~\ref{p=compmeas}, we establish a related
result, which will also be used in the proof of Theorem~\ref{t=main}.

\begin{proposition}\label{p=compmeas2} 
Let $\F$ be absolutely continuous with bounded Jacobians, and let $U$ be
a foliation box for $\F$ with smooth transversal $\tau$.  
Suppose that $\{Y_n\}_{n\geq 0}$ is a sequence
of measurable sets in $U$ with $c$-uniform fibers, for some $c\geq 1$.  
Then, for every $\F$-saturated measurable
set $X$, we have the equivalence:
$$
\lim_{n\to\infty} m(X:Y_n) = 1\,\Longleftrightarrow\, 
\lim_{n\to\infty}m_\tau(\tau_X:\tau_{Y_n}) = 1.
$$
\end{proposition}

\remark The hypothesis that $X$ is $\F$-saturated can be weakened: 
it suffices for $X \cap U$ to be a union of local leaves of $\F$ in $U$.

\bigskip

\noindent{\bf Proof of Proposition~\ref{p=compmeas2}.}
Let $X^*$ be the complement of $X$ in $M$. Then $X^*$ is also $\F$-saturated. The
proposition can be reformulated in terms
of $X^*$. We have to prove the equivalence:
$$
\lim_{n\to\infty} m(X^*:Y_n) = 0\,\Longleftrightarrow\, 
\lim_{n\to\infty}m_\tau(\tau_{X^*}:\tau_{Y_n}) = 0.
$$

For each $n$, let 
$$m_n = \inf_{q\in \tau_{Y_n}} m_\F(Y_n(q)).$$
Since the fibers of $Y_n$ are $c$-uniform, it
follows that $m_n>0$ for all $n$, and:
$$ m_n \leq {m_\F({Y_n}(q))} \leq cm_n,$$
for all $q\in \tau_{Y_n}$.  Absolute continuity implies that there
exists a $C\geq 1$ such that
$$
C^{-1} m(Y_n) \leq \int_{\tau} m_{\F}({Y_n}(q))\, dm_\tau(q)\leq C m(Y_n).
$$
Together, these inequalities imply that
\begin{eqnarray}\label{e=yn}
C^{-1}m_n m_\tau(\tau_{Y_n}) \leq  m(Y_n) \leq Cc m_n m_\tau(\tau_{Y_n}).
\end{eqnarray}

Since $X^*$ is
$\F$-saturated, the $\F$-fiber of $X^*\cap Y_n$ over a point 
$q\in \tau_{Y_n}$ is either empty or equal to $Y_n(q)$.  
Thus $X^*\cap Y_n$ also has
$c$-uniform fibers, and, as above, we obtain:
\begin{eqnarray}\label{e=syn}
C^{-1}m_n m_\tau(\tau_{X^*\cap Y_n}) \leq  m({X^\ast \cap Y_n}) 
\leq C c m_n m_\tau(\tau_{X^\ast\cap Y_n}).
\end{eqnarray}

Noting that $\tau_{X^\ast \cap Y_n} = \tau_{X^\ast} \cap \tau_{Y_n}$
and dividing the inequalities in (\ref{e=syn}) by those in 
(\ref{e=yn}), we obtain:
$$(C^2c)^{-1} m_\tau(\tau_{X^*}: \tau_{Y_n}) \leq m({X^*}:{Y_n}) 
\leq C^2c m_\tau(\tau_{X^*}: \tau_{Y_n}).$$
The result follows easily from this.
 \endproof

\bigskip
\noindent{\bf Proof of Proposition~\ref{p=compmeas}.
} Let $X'$ be an essential $\F$-saturate of $X$.  
Using Proposition~\ref{p=compmeas}, we have the equivalences:
\begin{eqnarray*}
\lim_{n\to\infty} m(X:Y_n) = 1 &\Longleftrightarrow& \lim_{n\to\infty} m(X':Y_n) = 1\\
 &\Longleftrightarrow&\lim_{n\to\infty}m_\tau(\tau_{X'}:\tau_{Y_n}) = 1\\
&\Longleftrightarrow&\lim_{n\to\infty}m_\tau(\tau_{X'}:\tau_{Z_n}) = 1\\
&\Longleftrightarrow& \lim_{n\to\infty} m(X':Z_n) = 1\\
&\Longleftrightarrow& \lim_{n\to\infty} m(X:Z_n) = 1.\\
\end{eqnarray*}
\endproof

\subsection{Fake invariant foliations}\label{s=fake}

The $su$-density sequences that we use in this proof will be constructed using
dynamical foliations with uniform continuity properties. If
$f$ happens to be
dynamically coherent, then we are free to use the foliations 
$\W^s,  \W^u, \W^{c}, \W^{cs}$, and
$\W^{cu}$ for these constructions.  Since we are not 
assuming dynamical coherence, we must find  substitutes for $\W^c, \W^{cs}$,
and $\W^{cu}$ to make our proof work in general.
It turns out to be simplest to find substitutes
for {\em all} invariant foliations $\W^s, \W^u, \W^{c}, \W^{cs}$, and
$\W^{cu}$. We call these substitutes ``fake invariant
foliations.''  
There are a few key places in the argument where we will have
to use the real invariant foliations $\W^u$ and $\W^s$, rather
than their fake counterparts.
We will indicate where this is the case.  The reader should
recall the choice of the
constant $R$ from Section~\ref{s=adapted}.

\begin{proposition}\label{p=localfol} Let $f:M\to M$ be a $C^1$
partially hyperbolic diffeomorphism. For any $\eps > 0$, there exist
cconstants $r$ and $r_1$ with  $R> r > r_1 >0$ such that, for every $p\in M$, 
the neighborhood $B(p,r)$ is foliated by
foliations $\hW^u_p$, $\hW^s_p$, $\hW^c_p$, $\hW^{cu}_p$ and $\hW^{cs}_p$ with the
following properties:
\begin{enumerate}
\item {\bf Almost tangency to invariant distributions:}  
For each $q\in B(p, r)$, and every
$\beta \in \{u,s, c, cu, cs\}$, the leaf $\hW^\beta_p(q)$ is $C^1$ and
the tangent space $T_q\hW^\beta_p(q)$ lies in a cone of radius
$\eps$ about $E^\beta(q)$.
\item {\bf Local invariance:} for $q\in B(p, r_1)$,
we have:
$$f(\hW^\beta_p(q,r_1)) \subset 
\hW^\beta_{f(p)}(f(q)),\,\hbox{ and }f^{-1}(\hW^\beta_p(q,r_1)) \subset 
\hW^\beta_{f^{-1}(p)}(f^{-1}(q)).$$
\item {\bf Exponential growth bounds at local scales:} The following hold for all $n\geq 0$. 
\begin{enumerate}
\item Suppose that $q_j\in B(p_j,r_1)$ for $0 \leq j \leq n-1$.  

If
$q' \in \hW^{s}_{p}(q,r_1)$,
then $q_n' \in \hW^{s}_{p}(q_n,r_1)$, and
$$
d(q_n,q_n') \leq \nu_n(p) d(q,q').
$$
If $q_j' \in \hW^{cs}_{p}(q_j,r_1)$ for $0 \leq j \leq n-1$, then 
$q_n' \in \hW^{cs}_{p}(q)$, and
 $$d(q_n,q_n') \leq \hat\gamma_n(p)^{-1} d(q,q').$$

\item Suppose that $q_{-j}\in B(p_{-j},r_1)$ for $0 \leq j \leq n-1$. 

If $q' \in \hW^{u}_{p}(q,r_1)$,
then $q_{-n}' \in \hW^{u}_{p}(q_{-n},r_1)$, and
$$d(q_{-n},q_{-n}') \leq \hat\nu_{-n}(p)^{-1} d(q,q').$$
If $q_{-j}' \in \hW^{cu}_{p}(q_{-j},r_1)$ for $0 \leq j \leq n-1$, then
 $q_{-n}' \in \hW^{cu}_{p}(q_{-n})$, and
$$d(q_{-n},q_{-n}') \leq \gamma_{-n}(p) d(q,q').$$
\end{enumerate}

\item {\bf Coherence:}  $\hW^s_p$ and $\hW^c_p$ subfoliate  $\hW^{cs}_p$; 
 $\hW^u_p$ and $\hW^c_p$ subfoliate  $\hW^{cu}_p$.

\item {\bf Uniqueness:}
$\hW^s_p(p)= \W^s(p,r)$, and $\hW^u_p(p)= \W^u(p,r)$.

\item {\bf Regularity:} If $f$ is $C^{1+\delta}$, then the foliations $\hW^u_p$, $\hW^s_p$, $\hW^c_p$, $\hW^{cu}_p$ and $\hW^{cs}_p$ and their tangent distributions are uniformly
H\"older continuous.

\item{\bf Regularity of the strong foliation inside weak leaves:}
If $f$ is $C^2$ and center bunched, then each leaf of $\hW^{cs}_p$ is $C^1$ foliated by leaves of the foliation $\hW^s_p$, and each leaf of $\hW^{cu}_p$ is 
$C^1$ foliated by leaves of the foliation $\hW^u_p$.  If $C^{1+\delta}$ and strongly center bunched, then the same conclusion holds.
\end{enumerate}
The regularity statements in 6. and 7. hold uniformly in $p\in M$.

\end{proposition}

\proof Suppose that $f$ is $C^k$, for some $k\geq 1$. After possibly reducing $\eps$, we can assume that inequalities
     (\ref{e=defn})--(6)  %\ref{e=bunch})
hold for unit vectors in the $\eps$-cones around the spaces in the partially hyperbolic splitting.

The construction of the leaves of $\hW^{cu}_p$ and $\hW^{cs}_p$ through
$p$ is essentially the same as the proof of the existence of pseudo-hyperbolic plaque families in \cite{HPS}.  They are obtained as fixed points of
graph transforms of a map that coincides with $f$ in a neighborhood of the
orbit of $p$.  We take the argument one step further and consider all fixed points
of these graph transforms in the entire neighborhood of $p$.

Our construction will be performed in two steps. In the first, we construct 
foliations of each tangent space $T_pM$. In the second step, we use the
exponential map $\exp_p$ to
project these foliations from a neighborhood of the origin in $T_pM$ to 
a neighborhood of $p$.

\medskip  

\noindent{\bf Step 1.} \
We choose an $r_0>0$ such that $\exp_p^{-1}$ is defined
on $B(p,2r_0)$.  For $r\in (0,r_0]$, we define,
in the standard way, a map:
\fig{
\object {a}  0,4 {TM};
\object {b}  8,4 {TM};
\object {c}  0,0 {M};
\object {d}  8,0  {M};
\arrow  {a}  {b}  {F_r};
\arrow  {c}  {d}  {f};
\arrow  {a}  {c}  {};
\arrow  {b}  {d}  {};}
\noindent which is uniformly $C^k$ on fibers, satisfying:
\begin{enumerate}
\item $F_r(p,v) = \exp_{f(p)}^{-1}\circ f \circ \exp_p(v)$, for $\|v\|\leq r$;                                                           
\item $F_r(p,v) = T_pf(v)$, for $\|v\| \geq 2r$; 
\item  
$\|F_r(p,\cdot) - T_p f(\cdot)\|_{C^1}\to 0$ as $r\to 0$, uniformly in $p$.
\end{enumerate}

Endowing $M$ with the discrete topology,
we regard $TM$ as the disjoint union of its fibers. 
Property 3. implies that, if $r$ is small enough, then $F_r$
is partially hyperbolic, and each bundle in the partially
hyperbolic splitting for $F_r$ at $v\in T_pM$ lies within
the $\eps/2$-cone about the corresponding subspace of $T_pM$
in the  partially
hyperbolic splitting for $f$
at $p$ (we are making the usual identification of $T_vT_pM$ with $T_pM$).
If $r$ is small enough, the equivalents  of inequalities 
  (3)--(6)
will hold for $TF_r$.

If $r$ is sufficiently small, 
standard graph transform arguments give stable, unstable, center-stable,
and center-unstable foliations for
$F_r$ inside
each $T_pM$. We obtain a center foliation by intersecting the leaves of the
center-unstable and center-stable foliations.  While $TM$ is not compact, all of the relevant estimates for $F_r$ are uniform, and it is this,
not compactness, that counts.  The proof of the Hadamard-Perron Theorem in \cite{KH} contains many of details
of this argument.

These foliations are unique because they arise
as fixed points of globally-defined contracting 
graph transform maps. Consequently,
the stable foliation subfoliates the center-stable, and the 
unstable subfoliates the center-unstable. 

We now discuss the regularity properties of these foliations
of $TM$.  Recall the standard 
method for determining the regularity of 
invariant bundles and foliations. Suppose that 
$TX = E_1 \oplus E_2$ is a $Tg$-invariant
splitting of the tangent bundle for a $C^k$ 
diffeomorphism $g:X\to X$ satisfying, for every $p\in X$, and every
unit vector $v\in T_pX$:
\begin{alignat}{4}\label{e=defn6}
\alpha_1(p) < \|&Tg v\|  < \beta_1(p),&\qquad \hbox{if } v\in E_1(p), \\
\alpha_2(p) < \|&Tg v\|   &\qquad \hbox{if } v\in E_2(p), 
\end{alignat}
where $0<\alpha_1(p) < \beta_1(p)  < \alpha_2(p)$.  Then the bundle
$E_2$ is $C^a$, for any $a\leq k-1$ that satisfies:
$$ \sup_{p\in M} \frac{\beta_1(p)}{\alpha_1(p)^a\alpha_2(p)} < 1.$$

When the functions $\alpha_1, \alpha_2, \beta_2$ are constant, 
this fact is classical --- see, e.g. the $C^r$ Section Theorem in \cite{HPS}.
The general case of this result appears more recently in the literature
\cite{ss, hass, amiethesis, PSW}.  This result extends, at least in part,
to give regularity of invariant foliations.  In particular, 
when there is a foliation $\F_2$ tangent to $E_2$ that arises at the unique 
fixed point of a nonlinear graph transform, then $\F_2$
is a $C^a$ foliation \cite{PSW, PSWc}. 
These results are proved in the compact case, but compactness is used only to obtain uniform estimates on the 
functions $\alpha_1,\beta_1,\alpha_2$ and the derivative of~$g$; the results carry over as long as such uniform estimates hold.

Our foliations of $TM$ have been constructed as 
the unique fixed points of graph transform maps. 
We can apply the above results to the $F_r$-invariant splittings of $TTM$ as 
the sum of the stable and  center-unstable bundles for $F_r$
and as the sum of the center-stable and unstable bundles for $F_r$. 
It follows immediately that both the center-unstable and unstable bundles and 
the corresponding foliations are H\"older continuous as long as  
$F_r$ is $C^{1+\delta}$ for some $\delta >0$. 
We obtain the H\"older continuity of the 
center-stable and stable bundles for $F_r$ and the corresponding foliations
by thinking of the same splittings as $F_r^{-1}$-invariant.
H\"older regularity of the center bundle and foliation is obtained by
noticing the the center is the intersection of the center-stable and center-unstable.

When $k\geq 2$, a similar estimate gives
the $C^1$ regularity of the unstable bundle along the leaves of
the center-unstable foliation.  The manifold $X$ is the disjoint
union of the leaves of the center-unstable foliation for $F_r$, 
$E_2$ is the unstable bundle, and $E_1$ is the center bundle.
We have: 
$$\alpha_1 = \gamma,\quad \beta_1 = \hat\gamma^{-1},\quad\hbox{ and }
\quad\alpha_2 = \hat\nu^{-1}.$$
The center bunching hypothesis $\hat\nu < \gamma\hat\gamma$ implies that
$$ \sup_{p\in M} \frac{\hat\nu(p)}{\gamma(p)^a\, \hat\gamma(p)} < 1,$$
for some $a>1$.
It follows that
$E_2 = E^u$ is a $C^1$ bundle over $X$, the leaves of the
center-unstable foliation.  Similarly, this argument shows
that the  center bunching hypothesis $\nu < \gamma\hat\gamma$ implies that
$E^s$ is a $C^1$ bundle over the leaves of the
center-stable foliation. There is an additional difficulty in the proof of 
this estimate, which is that $X$ is not a $C^k$ manifold in general,
even when $F_r$ is $C^k$.  This difficulty is dealt with in \cite{PSW, PSWc}.

When $1<k<2$, this type of estimate does not work at all.  
A different argument, working with holonomy maps instead of bundles, is presented in  \cite{BWhold}. The main result there implies
that under the strong center bunching hypothesis on $f$, 
which carries over to $F_r$,
the unstable foliation $C^1$-subfoliates the center-unstable, and the  
stable foliation $C^1$-subfoliates 
the center-stable.

\medskip

\noindent{\bf Step 2.} \
We now have foliations of $T_pM$, for
each $p\in M$.  We obtain the foliations $\hW_p^u, \hW_p^c, \hW_p^s, \hW_p^{cu}$, and $\hW_p^{cs}$ by applying the exponential map $\exp_p$ to the corresponding foliations of $T_pM$ inside the ball around the origin of radius $r$.

If $r$ is sufficiently small, then
the distribution $E^\beta(q)$ 
lies within the angular $\eps/2$-cone about the parallel translate of
$E^\beta(p)$, for every $\beta \in \{u,s, c, cu, cs\}$ and all $p, q$
with  $d(p,q)\leq r$.  Combining this fact with the preceding discussion, we obtain that property 1. holds if $r$ is sufficiently small.

Property 2. --- local invariance --- follows from invariance under $F_r$ of the foliations of $TM$ and the fact that $\exp_{f(p)}(F_r(p,v)) = f(\exp_p(p,v))$ provided $\|v\| \leq r$.

Having chosen $r$, we now choose $r_1$ small enough so that 
$f(B(p,2r_1)) \subset B(f(p),r)$ and 
$f^{-1}(B(p,2r_1)) \subset B(f^{-1}(p),r)$, and so that, for all
$q\in B(p,r_1)$,
\begin{eqnarray*}
q' \in \hW_p^s(q,r_1)&\Longrightarrow& d(f(q),f(q')) \leq \nu(p)\, d(q,q'),\\
q' \in \hW_p^u(q,r_1)&\Longrightarrow&d(f^{-1}(q),f^{-1}(q')) \leq \hat\nu(f^{-1}(p))\, d(q,q'),\\
q' \in \hW_p^{cs}(q,r_1)&\Longrightarrow& d(f(q),f(q')) \leq \hat\gamma(p)^{-1} \,d(q,q'),\qquad and\\
q' \in \hW^{cu}_{p}(q,r_1)&\Longrightarrow& d(f^{-1}(q),f^{-1}(q')) \leq \gamma(f^{-1}(p))^{-1}\, d(q,q').
\end{eqnarray*}

Property 3. --- exponential growth bounds at local scales --- is now
proved by an inductive argument similar to the proof of 
Lemma~\ref{l=distest}.  

Properties 4.-- 7. --- coherence, uniqueness, regularity and regularity of the strong foliation inside weak leaves --- follow immediately from the corresponding properties of the foliations of $TM$ discussed above.
\endproof

\remark Note that the system of local foliations constructed in Proposition~\ref{p=localfol} is not unique; it depends on the extension of $F_r$ outside of a neighborhood of the zero-section of $TM$.  Also note that, even when $f$ is
dynamically coherent, in general there is no reason to expect
the fake invariant foliations $\hW_p^{cs}, \hW_p^{cu}$, and $\hW_p^c$ 
to coincide with
the local leaves of the
real invariant foliations  $\W^{cs}, \W_p^{cu}$, and $\W_p^c$, even at $p$.
For the leaf $\hW^{cs}_p(p)$ to coincide with the local leaf of
$\W^{cs}$ through $p$, it is necessary that every iterate $f^{-n}(\W^{cs}(p_n,r)$ overflow the neighborhood $B(p,r)$.

For the rest of the paper, $\hW^s_p,\hW^c_p$, $\hW^s_p$, $\hW^{cs}_p$
and $\hW^{cu}_p$ will denote fake invariant foliations given by Proposition~\ref{p=localfol}, with $\eps>0$ much less than the angle between any two
of the subspaces in the partially hyperbolic splitting.
We may rescale the metric so that the radius $r_1$ in
conclusion 2. of Proposition~\ref{p=localfol} is much bigger than $1$. 
This will ensure that all of the
objects used in the rest of the paper are well-defined.  
We may assume that if $\max\{d(x,p),d(y,p)\} \leq 3$, then
$\hW^{cs}_p(x) \cap \hW^u_{p}(y)$,
$\hW^{cs}_p(x) \cap \W^u_{loc}(y)$, $\hW^{cu}_{p}(x) \cap \hW^s_{p}(y)$
and  $\hW^{cu}_{p}(x) \cap \W^s_{loc}(y)$
are single points. We denote by $\hm_a$ the measure $m_{\hW^a}$.

\subsection{A simple distortion lemma}\label{s=dist}

The next lemma will be used to compare values of H\"older cocycles at 
nearby points.

\begin{lemma}\label{l=dist} Let $\alpha:M\to {\bf R}$ be a positive H\"older continuous function, with exponent $\theta>0$.  Then there exists a constant $H>0$ such that the following holds, for all $p,q\in M$, $B>0$ and $n\geq 1$:
$$\sum_{i=0}^{n-1} d(p_i, q_i)^\theta \leq B\quad\Longrightarrow\quad e^{-{H}B} \leq \frac{\alpha_n(p)}{\alpha_n(q)} \leq e^{{H}B},$$
and
$$\sum_{i=1}^{n} d(p_{-i}, q_{-i})^\theta \leq B \quad\Longrightarrow\quad e^{-{H}B} \leq \frac{\alpha_{-n}(p)}{\alpha_{-n}(q)} \leq e^{{H}B}.$$
\end{lemma}

\proof
We prove the first part of the proposition.  The second part is proved similarly. The function $\log\alpha$ is also H\"older continuous with exponent $\theta$.
Let $H>0$ be the H\"older constant of $\log\alpha$, so that for all
$x,y\in M$:
$$\left|\log\alpha(x) - \log\alpha(y)\right| \leq H d(x,y)^\theta.$$
The desired inequalities are equivalent to:
\begin{eqnarray*}
\left|\log\alpha_n(p)-\log\alpha_n(q)\right| \leq HB.
\end{eqnarray*}
Expanding $\log\alpha_n$ as a series, we obtain:
\begin{eqnarray*}
\left|\log\alpha_n(p)-\log\alpha_n(q)\right|& \leq & 
\sum_{i=0}^{n-1} \left| \log\alpha(p_i) - \log\alpha(q_i)\right|\\
&\leq& H \sum_{i=0}^{n-1}  d(p_i,q_i)^\theta.\\
&\leq& HB,
\end{eqnarray*}
since
$$\sum_{i=0}^{n-1} d(p_i, q_i)^\theta \leq B,$$
by the hypothesis of the Proposition.
\endproof

\subsection{Thin neighborhoods of $\W^s(p,1)$}\label{s=thin}

We next identify, for each $n\geq 0$ and $p\in M$, 
a neighborhood of $p$ whose first $n$ iterates 
remain in a uniform neighborhood of the corresponding iterates
of $p$. We give an exponential estimate of the
size of the first $n$ iterates of the $n$th such neighborhood.
In our proof of Theorem~\ref{t=main} we construct sequences of 
geometric objects; the  $n$th term in the sequence of objects for any
$x\in \W^s(p,1)$ will
lie in the $n$th neighborhood of $p$.

Let $\sigma<1$ be a continuous function.  For $n\geq 0$ and $p\in M$, define the
set $S_{n,\sigma}(p)$ by:
$$S_{n,\sigma}(p) = \bigcup_{x\in \W^s(p,1)} \hW^c_p(x,\sigma_n(p)).$$

\begin{lemma}\label{l=sliceest}Suppose that $\sigma$ satisfies $\sigma < \min\{\hat\gamma, 1\}$.
Then
$$f^j(S_{n,\sigma}(p)) \subset B(p_j, 2),$$
for $j=0,\ldots, n$.

Further, there exist positive constants
$\kappa < 1$ and $C>0$ such that, for every $n\geq 0$,
$$f^j(S_{n,\sigma}(p)) \subset B(p_j, C\kappa^j),$$
for $j=0,\ldots, n$. 
\end{lemma}

\proof 
Suppose that $x\in \W^s(p,1)$
and  $y\in \hW^c_p(x,\sigma_n(p))$.  By part 3(a) of Proposition~\ref{p=localfol}, we then have
$$y_j\in \hW^c_{p_j}(x_j, \hat\gamma_{j}^{-1}(p)\sigma_{n}(p))
\subset \hW^c_{p_j}(x_j, 1) \subset B(p_j,2),$$
for $0\leq j\leq n$.  In fact, since
$\sigma < \min\{\hat\gamma, 1\}$,
the quantity $\hat\gamma_{j}^{-1}(p)\sigma_{n}(p) <
\hat\gamma_{j}^{-1}(p)\sigma_{j}(p)$ is exponentially small in $j$,
as is the diameter of $f^j(\W^s(p,1))$. This implies the second conclusion.
\endproof

Now let $\tau\leq 1$ be another continuous function.  For every $x\in S_{n,\sigma}(p)$, we have that $B(x_n, \tau_n(p)) \subset B(p_n,r)$, and
so the set 
$$T_{n,\sigma,\tau}(p) = f^{-n}\left(\bigcup_{z\in f^n(S_{n,\sigma}(p))}\hW^u_{p_n}(z,\tau_n(p))\cup\W^u(z,\tau_n(p))\right)$$
is well-defined.  Proposition~\ref{p=localfol} and Lemma~\ref{l=distest} imply
that the leaves of $\hW^u_{p_j}$ and $\W^u_{loc}$ are uniformly contracted
by $f^{-1}$ as long as they stay near the orbit of $p$; 
combining these facts with Lemma~\ref{l=sliceest}, we get:

\begin{lemma}\label{l=kappa} For every continuous function $\sigma$ satisfying
$\sigma < \min\{\hat\gamma, 1\}$, the set $T_{n,\sigma, 1}(p)$
satisfies
$$f^j(T_{n,\sigma,1}(p)) \subset B(p_j, 3),$$
for $j=0,\ldots, n$.

Further, for every such $\sigma$ and every
continuous function $\tau<1$, there exist positive constants
$\kappa < 1$ and $C>0$ such that, for every $n\geq 0$,
$$f^j(T_{n,\sigma,\tau}(p)) \subset B(p_j, C\kappa^j),$$
for $j=0,\ldots, n$.
\end{lemma}

\begin{figure}\label{f=tn}
\begin{center}
\psfrag{t1}{$T_{1,\sigma,\tau}(p)$}
\psfrag{tn}{$T_{n,\sigma,\tau}(p)$}
\psfrag{tn+1}{$T_{n+1,\sigma,\tau}(p)$}
\psfrag{ftn}{$f\left(T_{n,\sigma,\tau}(p)\right)$}
\psfrag{ft1}{$f\left(T_{1,\sigma,\tau}(p)\right)$}
\psfrag{ftn+1}{$f\left(T_{n+1,\sigma,\tau}(p)\right)$}
\psfrag{fntn}{$f^n\left(T_{n,\sigma,\tau}(p)\right)$}
\psfrag{fntn+1}{$f^n\left(T_{n+1,\sigma,\tau}(p)\right)$}
\psfrag{fn+1tn+1}{$f^{n+1}\left(T_{n+1,\sigma,\tau}(p)\right)$}
\psfrag{p}{$p$}
\psfrag{ws}{$\hW^s_p$}
\psfrag{wu}{$\hW^u_p$}
\psfrag{wc}{$\hW^c_p$}
\includegraphics{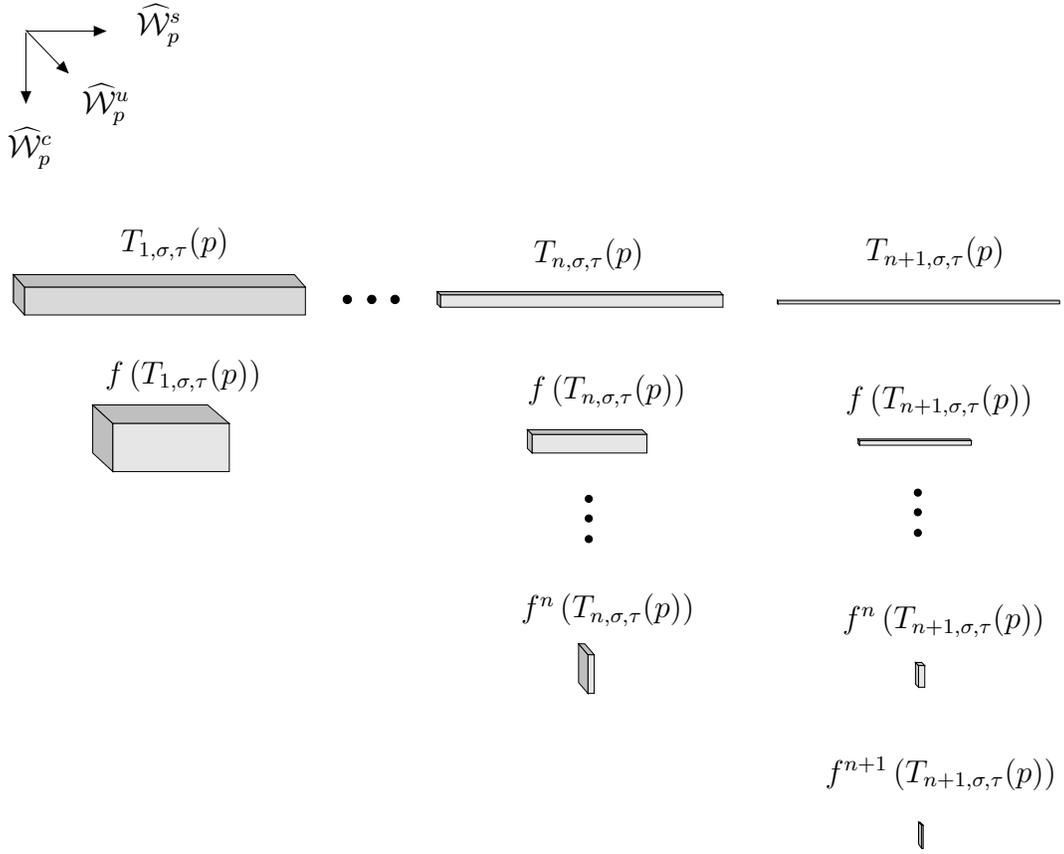}
\end{center}
\caption{Dimensions of the neighborhoods $T_{n,\sigma,\tau}(p)$ and their iterates}
\end{figure}

The dimensions of the neighborhoods $T_{n,\sigma,\tau}(p)$
and their iterates are illustrated in Figure~\ref{f=tn}, in the case where
$\sigma, \tau < 1$.

As simple corollary of this lemma
and Lemma~\ref{l=dist}, we then obtain:

\begin{lemma}\label{l=constalpha} Let $\alpha:M\to {\bf R}$ be a positive, uniformly H\"older continuous function, and let $\sigma, \tau$ be continuous
functions satisfying
$$\sigma < \min\{\hat\gamma, 1\}, \quad \tau < 1.$$  
Then there is a constant $C\geq 1$ such that,
for all $n\geq0$ and all $x, y\in T_{n,\sigma,\tau}(p)$, 
$$C^{-1}\leq \frac{\alpha_n(y)}{\alpha_n(x)} \leq C.$$
\end{lemma}

\section{The main theorem} \label{s=mainthm}

The properties of accessibility and essential accessibility can be reformulated
using the notion of saturation.
Accessibility means that a set which is both $\W^u$-saturated and $\W^s$-saturated
must be either empty or all of $M$. Essential accessibility means that a measurable set 
which is both $\W^u$-saturated and $\W^s$-saturated must
have either $0$ or full measure.

The central result of this paper is:

\begin{theorem}\label{t=reform} Let $f$ be $C^2$, 
partially hyperbolic and center bunched (or $C^{1+\delta}$ and strongly center
bunched).
Let $A$ be a measurable 
set that is both essentially $\W^u$-saturated and
essentially $\W^s$-saturated.
Then the set of Lebesgue density points of $A$ is 
$\W^u$-saturated and $\W^s$-saturated.
\end{theorem} 

The central result of Pugh and Shub in \cite{PS} is a version of 
Theorem~\ref{t=reform}, which involves a different notion of density point 
(defined in \cite{PS})  and
a slightly different hypothesis: 

{\em For $a = u$ or $s$, if $A$ is essentially 
$\W^a$-saturated, then the set of julienne density points of $A$ is
$\W^a$-saturated.}

 In contrast, Theorem~\ref{t=reform} requires that $A$ be 
{\em both} essentially $\W^u$-saturated and essentially $\W^s$-saturated 
in order to conclude anything.

Theorem~\ref{t=main} follows easily from Theorem~\ref{t=reform} by a 
version of the Hopf argument and a result of Brin and Pesin (see Section 2
of \cite{BPSW} for more details).

{\bf {\medskip}{\noindent}Proof of Theorem~\ref{t=main}. } 
To prove that $f$ is ergodic, it suffices to show that the Birkhoff averages of
continuous functions are almost everywhere constant.
Let $\varphi$ be a continuous function, and let
$$\hat\varphi_s(p) = 
\limsup_{n\to\infty} \frac{1}{n} \sum_{i=0}^n \varphi(f^i(p))
 \quad\hbox{ and }
\hat\varphi_u (p) = 
\limsup_{n\to\infty} \frac{1}{n} \sum_{i=0}^n \varphi(f^{-i}(p))
$$
be the forward and backward Birkhoff averages of $\varphi$
under $f$.  The function $\hat\varphi_s$ is constant along
$\W^s$-leaves, and  $\hat\varphi_u$ is constant along
$\W^u$-leaves.
It follows that for any $a\in \R$, the sets
$$A_s(a) = \hat\varphi_s^{-1}(-\infty,a]
\quad\hbox{ and }\quad A_u(a) = \hat\varphi_u^{-1}(-\infty,a].$$
are $\W^s$-saturated and
$\W^u$-saturated, respectively.  

The Birkhoff Ergodic Theorem 
implies that $\hat\varphi_s = \hat\varphi_u$ almost everywhere.
Consequently $m(A_s(a)\,\Delta\, A_u(a))=0$, so that 
the set $A(a) = A_u(a) \cap A_s(a)$ has 
$A_u(a)$ as an essential $\W^u$-saturate
and $A_s(a)$ as an essential $\W^s$-saturate. Thus
$A(a)$ satisfies 
the hypotheses of Theorem~\ref{t=reform}: it is
both essentially $\W^s$-saturated and essentially $\W^u$-saturated.

It follows from 
Theorem~\ref{t=reform}
that the set $\widehat A(a)$ of Lebesgue density points of $A(a)$ is  both
$\W^u$-saturated and $\W^s$-saturated. Essential accessibility implies that 
$\widehat A(a)$ has $0$ or full measure. But $m(A(a) \,\Delta\, \widehat A(a)) = 0$, so $A(a)$
itself has $0$ or full measure.  Since $a$ was arbitrary, it follows that
$\hat\varphi_s$ and $\hat\varphi_u$ are almost everywhere constant, and so
$f$ is ergodic.

To prove that $f$ has the Kolmogorov property, it suffices to show that all sets
in the Pinsker subalgebra~$\P$ have $0$ or full measure. 
According to Proposition 5.1 of \cite{BP1}, if $f$
is partially hyperbolic, then any set in $P \in \P$ is both essentially 
$\W^u$-saturated
and essentially $\W^s$-saturated. It again follows from Theorem~\ref{t=reform}
and essential accessibility that $P$ has $0$ or full measure.  \endproof

In order to prove Theorem~\ref{t=reform} it suffices to show that
the set of Lebesgue density points
of $A$ is $\W^s$-saturated; applying this result with $f$ replaced by $f^{-1}$
then shows that the set of Lebesgue density points
of $A$ is also $\W^u$-saturated. More precisely, it suffices to show that, for any $p\in M$,
if $x,x'\in \W^s(p,1)$, and $x$ is a Lebesgue density point of $A$,
then so is $x'$.
                                                                                
Let
$$N = \bigsqcup_{j\geq 0} B(p_j, r)$$ 
be the disjoint union, over $j\geq 0$, of the balls
$B(p_j, r)$, where $r\gg 1$ is given by Proposition~\ref{p=localfol}.    
Everything we do we do in the rest of this paper takes place inside $N$, and
we drop the dependence on $p$ where it is not confusing.  

We let $\hW^u$ be the locally invariant foliation of $N$ whose restriction to 
$B(p_j,r)$ is $\hW^u_{p_j}$.  Similarly we define foliations
$\hW^c, \hW^s, \hW^{cu}$ and $\hW^{cs}$.  By Proposition~\ref{p=localfol},
all of these foliations are
{\em uniformly} H{\"o}lder continuous, 
$\hW^{u}$ uniformly $C^1$ subfoliates $\hW^{cu}$,
and $\hW^{s}$ uniformly $C^1$ subfoliates $\hW^{cs}$. The foliations $\W^u$
and $\W^s$ induce foliations of
$N$, which will again be denoted by $\W^u$ and $\W^s$.
Hence $\W^s(p)$ is used to denote the
local leaf
$\W^s(p,r)$. Note that $\hW^s(p_j) = \W^s(p_j)$ and $\hW^u(p_j) = \W^u(p_j)$
for all $j\geq 0$. Note also that, since we are not assuming that
$f$ is dynamically coherent, there are no foliations
$\W^c, \W^{cu}$, or $\W^{cs}$.

Most leaves of these foliations
(with the notable exception the leaves of $\W^s$ and $\hW^s$ 
passing through the orbit
of $p$) are invariant under only finitely many iterates of
$f$, until their orbits leave $N$.  
In our proof, the geometric objects 
in $N$ that we need to iterate $n$ times always start out
in the neighborhood $T_{n,\sigma,1}$ defined in Section~\ref{s=thin} 
and therefore, by Lemma~\ref{l=kappa}, remain for $n$ iterates 
inside of $N$. As long as their orbits remain inside of $N$, the locally
invariant
foliations  $\hW^u, \hW^c, \hW^s, \hW^{cu}$ and $\hW^{cs}$ used to
construct these geometric objects are 
nearly indistinguishable from their invariant
counterparts $\W^u, \W^c, \W^s, \W^{cu}$ and $\W^{cs}$ in
the dynamically coherent setting.

These locally invariant
foliations differ from true invariant foliations in one key respect.
The measurable set $A$ in Theorem~\ref{t=reform} is
essentially $\W^s$-saturated and essentially $\W^u$-saturated. A set that is (essentially) $\W^u$-saturated need not be (essentially)
$\hW^u$-saturated, 
and a  set that is (essentially) $\W^s$-saturated need not be 
(essentially) $\hW^s$-saturated (nor do any of the converse implications hold). While it is possible  
to prove a version
of Theorem~\ref{t=reform} in which
$\W^u$ and $\W^s$ are replaced by $\hW^u$ and $\hW^s$, such a theorem 
is not enough to prove ergodicity using a Hopf argument.
This is because the Hopf argument uses infinitely many iterates of $f$ and
is thus global in nature.
Therefore, wherever we use explicitly the fact 
that  our set $A$ is $\W^u$-saturated, in particular, in our arguments that
use Proposition~\ref{p=compmeas}, we must
switch from using $\hW^u$ to using 
$\W^u$, and similarly for stable foliations.

Following Pugh and Shub, we consider for each $x\in \W^s(p,1)$ 
a sequence of sets,
called center-unstable juliennes, that
lie in the fake center-unstable manifold $\hW^{cu}(x)$ and shrink exponentially
as $n \to \infty$ while becoming increasingly thin in the $\hW^u$-direction.
In the Pugh-Shub construction, dynamical coherence is assumed,
and the true invariant foliations
$\W^{cu}$, $\W^c$ and $\W^u$ are used; here we use their
fake counterparts.
While objects we work with, such as center-unstable 
juliennes, will depend on the fake invariant foliations
given by Proposition~\ref{p=localfol}, the final conclusion of
Theorem~\ref{t=reform} does not, 
since we always have $\hW^s(p) = \W^s(p)$.

Recall the center bunching assumptions (\ref{e=bunch}):
\begin{eqnarray*}
  \nu < \gamma\hat\gamma  \qquad\hbox{and}\qquad    \hat\nu < \gamma\hat\gamma.
\end{eqnarray*}
In what follows we will use only the first of these inequalities.  The second inequality is used to prove $\W^u$ saturation of density points.

We choose continuous functions $\tau$ and $\sigma$ such that 
 $$
\nu < \tau < \sigma\gamma \quad\hbox{\rm and}\quad  \sigma < \min\{\hat\gamma,1\}.
 $$
Note that these inequalities  also imply that
$$\tau\hat\nu < \sigma\gamma\hat\nu <\sigma\gamma\hat\gamma \leq  \sigma.$$
The choice of $\sigma$ and $\tau$ with the desired properties is possible because of the center bunching assumption. The reader should think of
$\tau$ as being just a little bigger than $\nu$ and $\sigma$ as just a little bit less 
than $\min\{\hat\gamma,1\}$. 
The reader might also choose to keep in mind
the case where the functions $\nu, \hat\nu, \gamma$, and $\hat\gamma$
are constants, and where
$\tau$ and $\sigma$ can be chosen to be constant.  In this case the cocycles
$\tau_n$ and $\sigma_n$ are just the constants $\tau^n$ and $\sigma^n$.

Using the notation defined in Section~\ref{s=thin}, let
$S_n = S_{n,\sigma}(p)$, and let $T_n = T_{n,\sigma,\tau}(p)$.
For the rest of the paper, except where we indicate otherwise, cocycles will
be evaluated at the point $p$.  We will also drop the dependence on
$p$ from the notation; thus, if $\alpha$ is a cocycle, then
$\alpha_n(p)$ will be abbreviated to $\alpha_n$.

The center-unstable juliennes $\hJ^{cu}_n(x)$ that we construct will
be contained in $T_n$. 
We now describe the construction of center-unstable juliennes.

Define, for all $x\in \W^s(p,1)$,$$
\hB^c_n(x) = \hW^c(x,\sigma_n).
$$
Note that 
$$S_n = \bigcup_{y\in \W^s(p,1)} \hB^c_{n}(x).$$
For $y\in S_n$, we may then define two types of {\em unstable juliennes}:
$$
\hJ^u_n(y) = f^{-n}(\hW^u(y_n,\tau_n))
$$
and 
$$
J^u_n(y) = f^{-n}(\W^u(y_n,\tau_n)).
$$
Observe that for all $y,\in S_n$, the sets
$\hJ^u_n(y)$ and $J^u_n(y)$ are contained in $T_n$.

For each $x\in  \W^s(p,1)$ and $n\geq 0$, we then
define the {\em center-unstable julienne} centered at $x$ of order $n$:
$$
\hJ^{cu}_n(x) = \bigcup_{q \in \hB^c_n(x)} \hJ^u_n(q).
$$
Note that, by their construction, the sets
$\hJ^{cu}_n(x)$ are contained in $T_n$, for all $n\geq 0$ and
$x\in \W^s(p,1)$.
Note also that $\hJ^{cu}_n(x)$ is contained in the smooth submanifold 
$\hW^{cu}(x)$, which carries the restricted Riemannian volume
$\hm_{cu} = m_{\hW^{cu}}$, and $\hJ^{cu}_n(x)$ has positive $\hm_{cu}$--measure. If 
$J^u$ were used in place of $\hJ^u$ in the definition of
center-unstable julienne, the resulting set would not necessarily
be contained in a $C^1$ transversal to $\W^s$.

\begin{figure}\label{f=jules}
\begin{center}
\psfrag{bcnp}{$B^c_n\,\,\hB^c_n$}
\psfrag{fnws}{$\qquad J^u_n\,\,\,\,\hJ^u_n$}
\psfrag{p}{$p$}
\psfrag{pughshub}{Center-unstable julienne $J^{cu}_n(p)$ in \cite{PS}}
\psfrag{wc}{$\W^c$}
\psfrag{ws}{$\W^u$}
\psfrag{jcsn}{$\hJ^{cu}_n(p)$ in this paper.}
\includegraphics{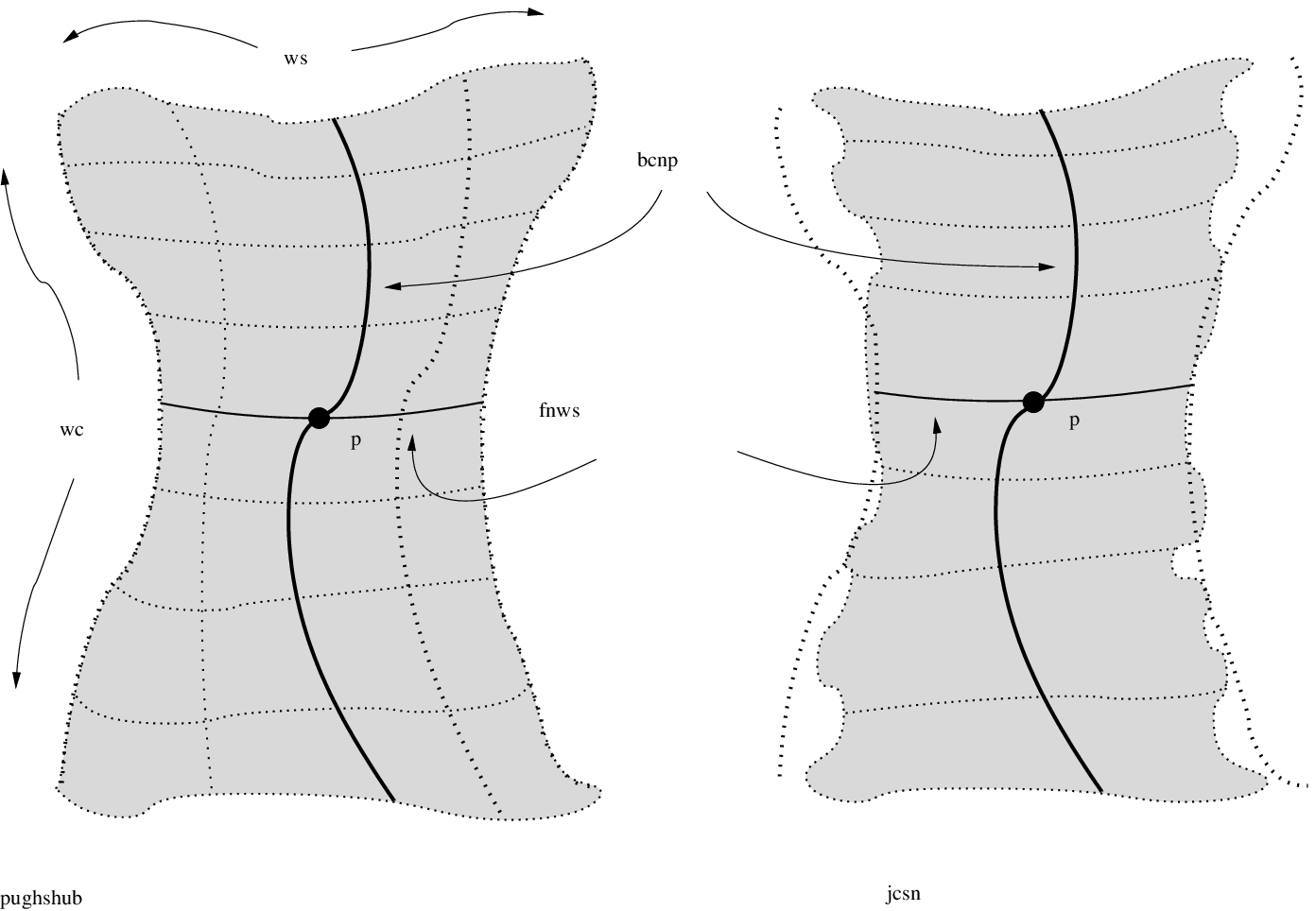}
\end{center}
\caption{Two types of center-unstable juliennes, when $\tau$ and $\sigma$ are constant.}
\end{figure}

Our $cu$-juliennes are closely related to, but not exactly the same as, those
of Pugh and Shub.
In the case where $\sigma$ and $\tau$ are constant functions, and
$f$ is dynamically coherent, their center-unstable julienne is 
the foliation product of
$\W^c(x,\sigma^n)$ and $f^{-n}(\W^u(x_{n}, \tau^n))$; see Figure~\ref{f=jules}.
In this case, the image under $f^{n}$ of our $J^{cu}_n(p)$ appears in \cite{PS}
as a tubelike approximation to the Pugh-Shub center-unstable postjulienne  of
rank~$n$. The results of \cite{PS} show that the  $cu$-juliennes defined
here and in \cite{PS} are comparable. Thus our $cu$-juliennes could be replaced by 
the Pugh-Shub $cu$-juliennes in Propositions~\ref{p=inclusion} and \ref{p=sliceden}.

As in \cite{PS}, the $cu$-juliennes have a quasi-conformality property: they 
are approximately preserved by holonomy along the stable foliation.

\begin{proposition}\label{p=inclusion} 
Let $x,x'\in \W^s(p,1)$, and let
$h^s: \hW^{cu}(x) \to \hW^{cu}(x')$  be the holonomy map
induced by the stable foliation
$\W^s$.
Then the sequences $h^s(\hJ^{cu}_n(x))$ and $\hJ^{cu}_{n}(x')$
are comparable.
\end{proposition}

There are estimates on the volumes of unstable and center-unstable juliennes.

\begin{proposition}\label{p=julmeas}  There exist $\delta > 0$ and
$c \geq 1$ such that, for all $x\in \W^s(p,1)$, and all $q,q'\in S_n$,
the following hold, for all $n\geq 0$:
$$
c^{-1}\leq \frac{\hm_{u}(\hJ^u_n(q))}{\hm_{u}(\hJ^u_n(q'))}\leq c,
$$
$$
c^{-1}\leq \frac{m_{u}(J^u_n(q))}{m_{u}(J^u_n(q'))}\leq c,
$$
$${\hm_{u}(\hJ^{u}_{n+1}(q))}\geq \delta {\hm_{u}(\hJ^{u}_{n}(q))},$$
and
$${\hm_{cu}(\hJ^{cu}_{n+1}(x))}\geq \delta {\hm_{cu}(\hJ^{cu}_{n}(x))}.$$
\end{proposition}

The final crucial property of the $cu$-juliennes is that, {\em
for the sets that
appear in the proof of Theorem~\ref{t=reform}},
Lebesgue density points are precisely $cu$-julienne density points.

\begin{proposition}\label{p=sliceden} 
Let $X$ be a measurable set that is both $\W^s$-saturated and
 essentially $\W^u$-saturated.  Then  $x\in \W^s(p)$ 
is a Lebesgue density point of $X$ if and only if:
$$
\lim_{n\to\infty} \hm_{cu}(X:\hJ^{cu}_n(x)) = 1.
$$
\end{proposition}

\remark Pugh and Shub show that
Lebesgue almost every point of {\em any} measurable set is
a $cu$-julienne density point.  In their argument they prove a Vitali covering
lemma for their juliennes.  This argument accounts for their
definition of $cu$-juliennes as a foliation product and for the stronger 
bunching hypothesis in their main result.  We do not know whether their result,
specifically  Theorem 7.1 of \cite{PS},
still holds under our weaker bunching hypothesis, or whether it holds at 
all in the absence of dynamical coherence (although it can
be shown that in the dynamically coherent, symmetrized setting that they
consider, 
their hypothesis can
be weakened from $\nu < \gamma^{2+2/\theta_0}$
to $\nu < \gamma^{1/\theta_0}$).

The proof of Proposition~\ref{p=inclusion} is essentially contained
in \cite{PS}. For completeness we reproduce the argument
in the next section.  Proposition~\ref{p=julmeas} is proved in Section~\ref{s=julmeas},
and the proof of Proposition~\ref{p=sliceden} is in the final section.  
We now use these three propositions to prove the main result.

{\bf {\medskip}{\noindent}Proof of Theorem~\ref{t=reform}. }
As we noted above, it suffices to show that the Lebesgue
density points of $A$ are $\W^s$-saturated;
to see that the Lebesgue density points of 
$A$ are $\W^u$-saturated, just consider $f^{-1}$ instead of~$f$.
Let $A^s$ be an essential $\W^s$-saturate of $A$.  
Since $m(A \,\Delta\, A^s) = 0$,
the Lebesgue density points of $A$ are precisely the same as those of $A^s$.
Fix $p\in M$ and suppose that $x\in \W^s(p,1)$ is a
Lebesgue density point of $A^s$.  Proposition~\ref{p=sliceden} 
implies that $x$ is a $cu$-julienne density point of $A^s$. 

To finish the proof, we show that every $x'\in W^s(p,1)$ is
a  $cu$-julienne density point of $A^s$.  Then 
by Proposition~\ref{p=sliceden}, every $x'\in W^s(p,1)$ is
a Lebesgue density point of $A^s$. The Lebesgue density
points of $A^s$, and hence of $A$, are therefore $\W^s$-saturated.  

Let $h^s: \hW^{cu}(x) \to \hW^{cu}(x')$  be the holonomy map
induced by the stable foliation $\W^s$.
The sequence $h^s(\hJ^{cu}_n(x))\subset \hW^{cu}(x')$ nests at $x'$.  
Transverse absolute continuity of $h^s$ with bounded 
Jacobians implies that
$$\lim_{n\to\infty} \hm_{cu}(A^s:\hJ^{cu}_n(x)) = 1\quad\Longleftrightarrow\quad
\lim_{n\to\infty} \hm_{cu}(h^s(A^s):h^s(\hJ^{cu}_n(x))) = 1.$$
Since $A^s$ is $s$-saturated, we then have:
$$\lim_{n\to\infty} \hm_{cu}(A^s:\hJ^{cu}_n(x)) = 1\quad\Longleftrightarrow\quad 
\lim_{n\to\infty} \hm_{cu}(A^s:h^s(\hJ^{cu}_n(x))) = 1.$$
Since we are assuming that $x$ is a $cu$-julienne density point of $A^s$, we thus have
$$\lim_{n\to\infty} \hm_{cu}(A^s:h^s(\hJ^{cu}_n(x))) = 1.$$

Working inside of $\hW^{cu}(x')$, we will apply Lemma~\ref{l=compreg}
to the sequences $h^s(\hJ^{cu}_n(x))$ and
$\hJ^{cu}_n(x')$, which both nest at $x'$.  
Proposition~\ref{p=inclusion} implies that these sequences are comparable.
Proposition~\ref{p=julmeas} implies that $\hJ^{cu}_n(x')$ is regular, 
with respect to the induced
Riemannian measure $\hm_{cu}$ on $\hW^{cu}(x')$.
Lemma~\ref{l=compreg} now tells us that 
$$\lim_{n\to\infty} \hm_{cu}(A^s:h^s(\hJ^{cu}_n(x))) = 1\quad\Longleftrightarrow\quad\lim_{n\to\infty} \hm_{cu}(A^s:\hJ^{cu}_n(x')) = 1,$$
and so $x'$ is a $cu$-julienne density point of $A^s$. It follows from Proposition~\ref{p=sliceden} that $x'$ is a Lebesgue density point of $A^s$, and thus of $A$. \endproof

\section{Julienne quasiconformality} \label{s=jqconf}
We adapt the proof of Theorem 4.4 in \cite{PS} to prove Proposition~\ref{p=inclusion}.
It will suffice to show that $k$ can be chosen so that
\begin{equation}\label{e=inc}
h^s(\hJ^{cu}_n(x))\subseteq \hJ^{cu}_{n-k}(x'),
\end{equation}
for all $n\geq k$, whenever $x$ and $x'$ satisfy the hypotheses of the
proposition.  The hypotheses of the proposition treat $x$ and $x'$ symmetrically, 
so we can then reverse their roles to obtain:
$$\overline{h}^s(\hJ^{cu}_n(x'))\subseteq \hJ^{cu}_{n-k}(x),$$
for all $n\geq k$, where $\overline{h}^s:\hW^{cu}_{loc}(x')\to \hW^{cu}(x)$
is the holonomy induced also by the stable foliation $\W^s$.  Since
 $\overline{h}^s$ and $h^s$ are inverses, we then obtain:
$$
\hJ^{cu}_n(x')\subseteq {h}^s(J^{cu}_{n-k}(x)),
$$
for all $n\geq k$.

In order to prove that $k$ can be chosen so that (\ref{e=inc}) holds,
we need two lemmas.

\begin{lemma}\label{l=centhol} There
exists a positive integer $k_1$ such that, for all $x, x'\in \W^s(p)$,
$$
\hat{h}^s(\hB^c_n(x)) \subseteq \hB^c_{n-k_1}(x'),
$$
for all $n\geq k_1$, where $\hat{h}^s:\hW^{cu}_{loc}(x)\to \hW^{cu}(x')$
is the local $\hW^s$ holonomy.
\end{lemma}

%\noindent{\bf Proof of Lemma~\ref{l=centhol}.} 
\proof Proposition~\ref{p=localfol} implies that $\hat{h}^s$ is 
$L$-Lipschitz, for some $L\geq 1$. Therefore the
image of $\hW^c(x,\sigma_n)$ under $\hat{h}^s$
is contained in $\hW^c(x',L\sigma_n)\subseteq \hW^c(x',\sigma_{n-k_1})$,
 for any
$k_1$ large enough so that $\sigma_{-k_1} > L$. \endproof

\begin{figure}[h]\label{f=holon}
\begin{center}
\psfrag{ws}{$\hW^u$}
\psfrag{hws}{$\hW^s$}
\psfrag{wu}{$\W^s$}
\psfrag{wc}{$\hW^c$}
\psfrag{p}{$p$}
\psfrag{hb}{$\hB^c_n(p)$}
\psfrag{x}{$q$}
\psfrag{y}{$y$}
\psfrag{xn}{$q_{n}$}
\psfrag{sn}{$S_{n}$}
\psfrag{yn}{$y_{n}$}
\psfrag{xprime}{$q'$}
\psfrag{yprime}{$y'$}
\psfrag{zprime}{$z'$}
\psfrag{xnprime}{$q'_{n}$}
\psfrag{ynprime}{$y'_{n}$}
\psfrag{znprime}{$z'_{n}$}
\psfrag{fminusn}{$f^{n}$}
\psfrag{Osigman}{$O(\tau_n)$}
\psfrag{Onuhatn}{$O(\nu_n)$}
\includegraphics{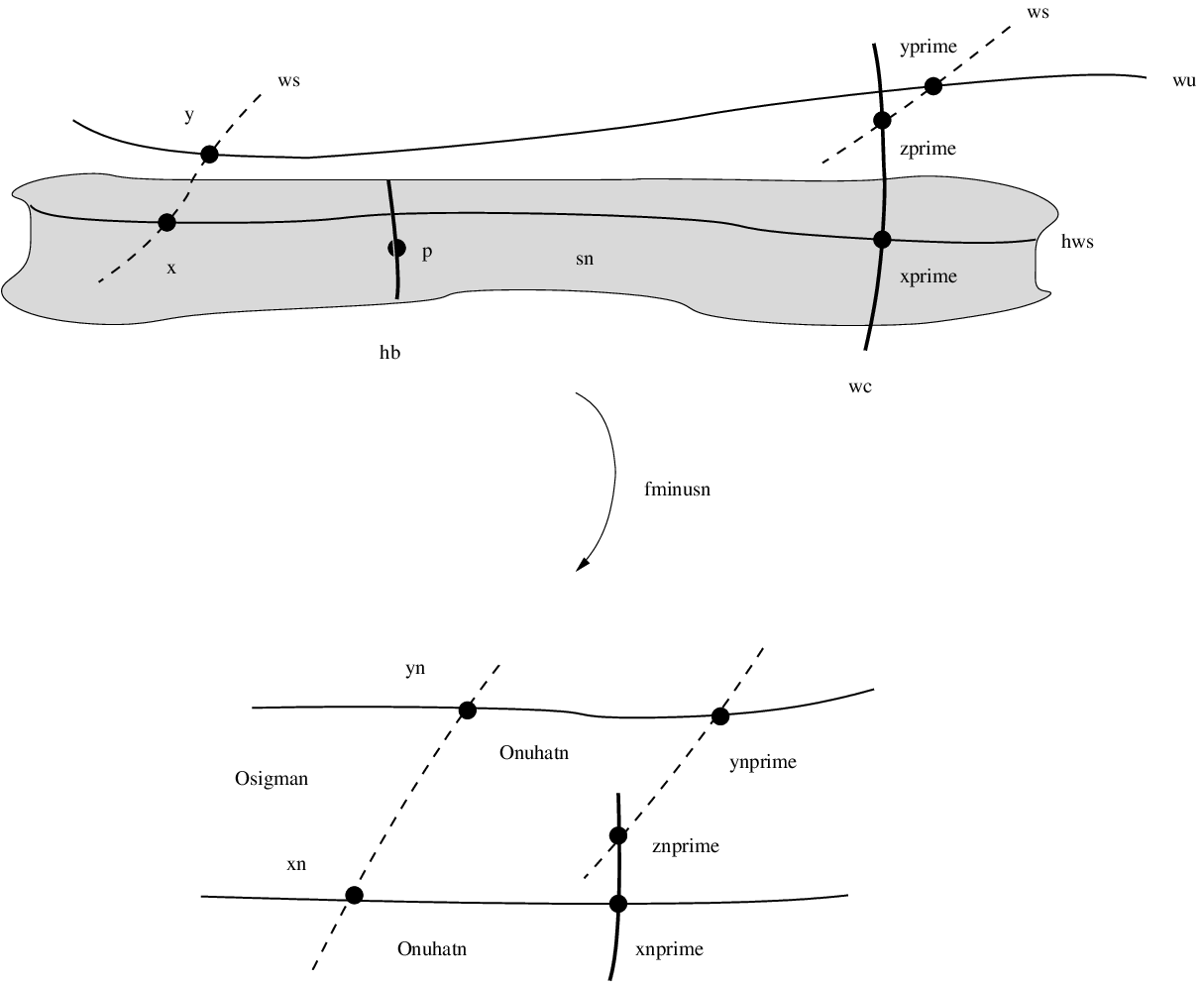}
\end{center}
\caption{Picture for the proof of Lemma~\ref{l=holon}.}
\end{figure}

\begin{lemma}\label{l=holon} There exists a positive integer 
$k_2$ such that the
following holds for every integer  
$n\geq k_2$. Suppose $q,q'\in S_n$ with $q' \in \hW^s(q)$.  
Let $y\in \hJ^u_n(q)$, and let $y'$ be the image of $y$ under $\W^s$ holonomy
from $\hW^{cu}_{loc}(q)$ to $\hW^{cu}(q')$.  Then
$$
y' \in \hJ^{u}_{n-k_2}(z'),
$$
for some $z'\in \hW^c(q', \sigma_{n-k_2})$.
\end{lemma}

\remark Note that two types of holonomy maps appear in Lemma~\ref{l=holon}:
the point $q'$ is the image of $q$ under $\hW^s$ holonomy between
$\hW^{cu}_{loc}(x)$ and $\hW^{cu}(x')$, whereas $y'$ is the image of 
$y$ under $\W^s$ holonomy. 

\bigskip

\noindent{\bf Proof of Lemma~\ref{l=holon}.}  
Let $z'$ be the unique point
in $\hW^u(y')\cap \hW^c(q')$.  It is not hard to see
that $z'_j\in N$, for $j=0,\ldots, n-1$ and that 
$z'_{n}$ is the unique point
in $\hW^u(y'_{n})\cap \hW^c(q'_{n})$. 
It will suffice to prove that $d(y_{n}',z_{n}') 
= O(\tau_n)$ and
$d(q',z') = O(\sigma_n)$.

We have $d(q_{n}, y_{n}) \leq \tau_n$ because 
$y\in f^{-n}(\W^u(q_{n},\tau_n))$. By Proposition~\ref{p=localfol}, 3(a),
we also have
that 
$d(q_{n},q'_{n}) = O(\nu_n)$ and  $d(y_{n},y'_{n}) = O(\nu_n)$, since
$d(q,q')$ and $d(y,y')$ are both $O(1)$. Note that
$q_n$ and $z_n'$ are, respectively, the images of $y_n$ and $y_n'$ 
under $\hW^u$-hononomy between $\hW^{cs}_{loc}(y_n)$ and 
$\hW^{cs}(q_n)$. Uniform continuity of the foliations
$\hW^u$ and $\hW^{cs}$ implies that
$$d(y_n',z_n') = O(\max\{d(q_n,y_n), d(y_n,y_n')\}) = O(\tau_n),$$
since $\nu < \tau$. 

We next show that $d(q',z') = O(\sigma_n)$.
By the triangle inequality,
$$d(q_{n}',z_{n}') \leq d(q_n',q_n) + d(q_n,y_n) + d(y_n, y_n') + d(y_n', z_n').$$
All four of the quantities on the right-hand side are easily seen to be 
$O(\tau_n)$.
Since $q_{n}'$ and $z_{n}'$ lie in the same $\hW^c$-leaf at 
distance $O(\tau_n)$, Proposition~\ref{p=localfol}
now implies that
$d(q',z') = O(\gamma_{-n}\tau_n)$. 
But $\tau$ and $\sigma$ were chosen so that $\tau < \gamma\sigma$.  
Hence $\gamma_{-n}\tau_n< \sigma_n$ and 
$d(q',z') = O(\sigma_n)$, as desired. \endproof

{\bf {\medskip}{\noindent}Proof of Proposition~\ref{p=inclusion}. }
As noted above, it suffices to prove the inclusion~(\ref{e=inc}).
For $q \in B^c_n(x)$, let $q' = \hat{h}^s(q)$. Then $q' \in B^c_{n-k_1}(p')$
by Lemma~\ref{l=centhol}. Hence $q,q'\in S_{n-k_1}$
 and we can apply Lemma~\ref{l=holon} to obtain
$$
h^s(J^{cu}_n(x)) \subseteq \bigcup_{z \in Q} J^{u}_{n-k_2}(z),
$$
where 
$$Q=\bigcup_{q' \in B^c_{n-k_1}(x')} B^{c}_{n-k_2}(q').$$
For $k\geq k_2$, we have:
$$\bigcup_{z \in Q} J^{u}_{n-k_2}(z)\subseteq 
\bigcup_{z \in Q} J^{u}_{n-k}(z).$$
It therefore suffices to find $k\geq k_2$ such that $Q\subseteq B^c_{n-k}(x')$.
This latter inclusion holds if:
$$
\sigma_{n-k_1} + \sigma_{n-k_2} \leq \sigma_{n-k},$$
which is obviously true for all $n\geq k$, if $k$ is sufficiently large.
 \endproof

\section{Julienne measure}\label{s=julmeas}
 We next prove Proposition~\ref{p=julmeas}. 
Continuity of $\hW^u$ implies that there exists $C_1\geq 1$ such that
\begin{eqnarray}\label{e=init}
C_1^{-1} \leq  
   \frac{\hm_{u}(\hW^u(q_{n},\tau_n))}{\hm_{u}(\hW^u(q_{n}',\tau_n))} 
             \leq C_1,
\end{eqnarray}
for all $q,q'\in S_n$.

Let $\hE^s, \hE^c$, and $\hE^u$ be the tangent distributions to the
leaves of $\hW^s, \hW^c$, and $\hW^u$, respectively. They are H\"older
continuous by Proposition~\ref{p=localfol}, part 6.   Furthermore, 
the restrictions of
these distributions to $T_n$ are invariant under $Tf^j$, for $j=1,\ldots n$.
We next observe that the Jacobian $\hbox{Jac}(Tf^{n}\vert_{\hE^u})$ 
is nearly constant when
restricted to the set $T_n$.
More precisely, we have:
\begin{lemma}\label{l=jacdist}
There exists $C_2\geq 1$ such that, for all $n\geq 1$, and
all $y,y'\in T_n$,
\begin{eqnarray*}
C_2^{-1}\leq \frac{\hbox{\rm Jac}(Tf^{n}\vert_{\hE^u})(y)}
                  {\hbox{\rm Jac}(Tf^{n}\vert_{\hE^u})(y')} \leq C_2.
\end{eqnarray*}
\end{lemma}

\proof
By the Chain Rule, these inequalities follow from Lemma~\ref{l=constalpha} with
$\alpha = \hbox{Jac}(Tf\vert_{\hE^u})$.\endproof

Let $q\in  S_n$, and let $X\subseteq \hJ^u_n(q)$ be a measurable
set (such as $\hJ^u_n(q)$ itself). Then:
\begin{eqnarray*}
\hm_u(f^n(X)) &=& \int_{X}
       \hbox{\rm Jac}(Tf^{n}\vert_{\hE^u})(x)\,d\hm_u(x).
\end{eqnarray*}
From this and Lemma~\ref{l=jacdist} we then obtain:
\begin{lemma}\label{l=generalcompare} There exists
$C_3>0$ such that, for all $n\geq 0$, for
any $q,q'\in S_n$, and any
measurable sets $X\subset \hJ^u_n(q), X'\subset \hJ^u_n(q')$,
we have:
$$C_3^{-1}\frac{\hm_u(f^n(X))}{\hm_u(f^{n}(X'))} \leq \frac{\hm_u(X)}{\hm_u(X')}
\leq C_3 \frac{\hm_u(f^n(X))}{\hm_u(f^{n}(X'))}.$$
\end{lemma}

Recall that $f^n(\hJ^u_n(q)) = \hW^u(q_n,\tau_n)$, for
$q\in S_n$.
The first conclusion of Proposition~\ref{p=julmeas} now
follows from (\ref{e=init}) and Lemma~\ref{l=generalcompare} with
$X= \hJ^u_n(q)$ and $X' = \hJ^u_n(q')$.

The second conclusion is proved similarly.

We next show that there exists $\delta>0$ such that
\begin{eqnarray}\label{e=juratio}
\frac{\hm_{u}(\hJ^{u}_{n+1}(q))}{\hm_{u}(\hJ^{u}_{n}(q))} \geq \delta,
\end{eqnarray}
for all $n\geq 0$ and all $q\in S_n$.
To obtain (\ref{e=juratio}), we will apply
Lemma~\ref{l=generalcompare} with $q=q'$, $X=\hJ^u_{n+1}(q)$,
and $X' = \hJ^u_n(q)$.  This gives us:
\begin{eqnarray*}
\frac{\hm_u(\hJ^{u}_{n+1}(q))}{\hm_u(\hJ^{u}_{n}(q))} & \geq&
C_3^{-1} \frac{\hm_u(f^{n}(\hJ^{u}_{n+1}(q)))}{\hm_u(f^{n}(\hJ^{u}_{n}(q)))}.
\end{eqnarray*}
But $f^{n}(\hJ^{u}_{n+1}(q)) =
f^{-1} (\hW^{u}(q_{n+1},\tau_{n+1}))$
and  $f^{n}(\hJ^{u}_{n}(q))=
\hW^u(q_n, \tau_n)$,
and hence:
\begin{eqnarray*}
\frac{\hm_u(f^{n}(\hJ^{u}_{n+1}(q)))}{\hm_u(f^{n}(\hJ^{u}_{n}(q)))} &= &
\frac{\hm_u(f^{-1} (\hW^{u}(q_{n+1},\tau_{n+1})))}{\hm_u(\hW^u(q_n, \tau_n))}.\\
\end{eqnarray*}
This ratio is uniformly bounded below away from $0$,
since $f^{-1}$ is a diffeomorphism, the leaves of $\hW^u$ are uniformly smooth,
 and the ratio
$\tau_{n+1}/\tau_{n} = \tau(p_n)$ is uniformly bounded below away from $0$.

To prove the final claim, we begin by observing that,
considered as a subset of $\hW^{cu}(x)$, the set $\hJ^{cu}_n(x)$ fibers over $\hB^c_{n}(x)$ with $\hW^u$-fibers $\hJ^u_n(q)$.  We have just proved that these fibers are $c$-uniform.
Since  $\sigma_{n+1}/\sigma_{n} = \sigma(p_n)$ is uniformly bounded away from $0$, the ratio
$$\frac{\hm_c(\hB^c_{n+1}(x))}{\hm_c(\hB^c_{n}(x))} = \frac{\hW^{c}(x,\sigma_{n+1})}{\hW^{c}(x,\sigma_{n})} $$
is bounded away from $0$, uniformly in $x$ and $n$.
Thus the sequence of bases $\hB^c_{n}(x)$ of $\hJ^{cu}_n(x)$ 
is regular in the induced Riemannian volume $\hm_c$.
Proposition~\ref{p=localfol}, part 7. 
implies that $\hW^u$ $C^1$ subfoliates
$\hW^{cu}$; in particular, considered as a subfoliation of $\hW^{cu}(x)$,
$\hW^u$ is absolutely continuous with bounded
Jacobians. Proposition~\ref{p=unifreg} implies
that the sequence $\hJ^{cu}_n(x)$ is regular, with respect to the 
induced Riemannian measure $\hm_{cu}$.  This proves the final claim.
\endproof

\section{Julienne density}\label{s=density}
We now come to the proof of Proposition~\ref{p=sliceden}. We must show that if
a measurable set $X$ is both $\W^s$-saturated and
 essentially $\W^u$-saturated, then a point 
 $x\in \W^s(p,1)$ is a Lebesgue density point of $X$ if and only if
$$
\lim_{n\to\infty} \hm_{cu}(X:\hJ^{cu}_n(x)) = 1.
$$

We will establish the following chain of equivalences:
\begin{eqnarray*}
 \hbox{$x$ is a Lebesgue density point of $X$} 
  &\Longleftrightarrow &
              \lim_{n\to\infty} m(X:B_n(x)) = 1 \\
  &\Longleftrightarrow &
              \lim_{n\to\infty} m(X:C_n(x)) = 1 \\
  &\Longleftrightarrow &
              \lim_{n\to\infty} m(X:D_n(x)) = 1 \\
  &\Longleftrightarrow &
              \lim_{n\to\infty} m(X:E_n(x)) = 1 \\
    &\Longleftrightarrow &
              \lim_{n\to\infty} m(X:F_n(x)) = 1 \\
  &\Longleftrightarrow &
              \lim_{n\to\infty} m(X:G_n(x)) = 1 \\
  &\Longleftrightarrow &
              \lim_{n\to\infty} \hm_{cu}(X:\hJ^{cu}_n(x)) = 1.
 \end{eqnarray*}
Before defining the sets $B_n(x)$ through $G_n(x)$, we outline the general 
scheme of the proof. After verifying the first equivalence, we prove that $B_n(x)$ is regular, and that
$B_n(x)$ and $C_n(x)$ are comparable.  The second equivalence then follows
from Lemma~\ref{l=compreg}. The sets $C_n(x)$ and $D_n(x)$ both fiber
over the same base in $\hW^{cs}$, with $c$-uniform fibers in $\W^u$,
so the third equivalence follows from Proposition~\ref{p=compmeas}. We prove that the sets
$D_n(x), E_n(x), F_n(x)$, and $G_n(x)$ are all comparable, and that
$G_n(x)$ is a regular sequence.  Equivalences 4-6 then follow from Lemma~\ref{l=compreg}.  Finally, $G_n(x)$ fibers over $\hJ^{cu}_n(x)$, with
$c$-uniform $\W^s$-fibers, and so the final equivalence follows from
Proposition~\ref{p=compmeas2}. This final step uses $\W^s$-saturation of $X$.

The sets $B_n(x)$ through $G_n(x)$ are defined as follows.
The set $B_n(x)$ is a Riemannian ball in the manifold $\hW^c(x)$:
$$
B_n(x) = B(x,\sigma_n).
$$
The sets  $C_n(x)$, $D_n(x)$ and $E_n(x)$ 
will fiber over the same base $D^{cs}_n(x)$, where
$$
 D^{cs}_n(x) = \bigcup_{x'\in \hW^s(x,\sigma_n)}\hB^c_n(x').
 $$ 
Proposition~\ref{p=localfol}, part 4. implies that
$D^{cs}(x)$ is contained in the $C^1$ submanifold $\hW^{cs}(x)$;
it is comparable to the disk $\hW^{cs}(x,\sigma_n)$.  
Let 
$$C_n(x)  
=\bigcup_{q\in D^{cs}_n(x)} \W^u(q,\sigma_n),$$
and let
$$
D_n(x) = \bigcup_{q\in D^{cs}_n(x)} J^u_n(q).$$
The set $E_n(x)$ is nearly identical to $D_n(x)$, with the crucial difference 
that the $J^u_n$-fibers are replaced with
$\hJ^u_n$-fibers:
$$
E_n(x) = \bigcup_{q\in D^{cs}_n(x)} \hJ^u_n(q) = 
\bigcup_{x' \in \hW^s(x,\sigma_n)} \hJ^{cu}_n(x')
= \bigcup_{x' \in \W^s(x,\sigma_n)} \hJ^{cu}_n(x').$$
The rightmost equality follows
from the fact that $\hW^s(x,\sigma_n) =\W^s(x,\sigma_n)$,
for all $x\in \W^s(p,1)$ (Proposition~\ref{p=localfol}, part 5.)

We define $F_n(x)$ to be the foliation product of $\hJ^{cu}_n(x)$ and $\W^s(x,\sigma_n)$:
$$
F_n(x) = \bigcup_{q\in \hJ^{cu}_n(x), \,\, q'\in \W^s(x,\sigma_n)} \W^s(q)\cap
\hW^{cu}(q').
$$
This definition makes sense since the 
foliations $\hW^{cu}$ and $\W^s$ are transverse.
  Finally, let
$$
G_n(x) = \bigcup_{q \in \hJ^{cu}_n(x)} \W^s(q,\sigma_n).$$
It is in the transition from $D_n$ to
 $E_n$ that the exchange between the measure-theoretically useful foliation $\W^u$ and the geometrically useful foliation $\hW^u$
takes place. The definition of $F_n$ is where we first use
the foliation $\W^s$.

\begin{figure}\label{f=efg}
 \begin{center}~
\psfrag{wu}{$\hW^s$}
\psfrag{wcs}{$\hW^{cu}$}
\psfrag{p}{$x$}
\psfrag{E}{$E_n(x)$}
\psfrag{F}{$F_n(x)$}
\psfrag{G}{$G_n(x)$}
\psfrag{jcsn}{$\hJ^{cu}_n(x)$}
\psfrag{wuprhon}{$\hW^s(x,\sigma_n)$}
\includegraphics{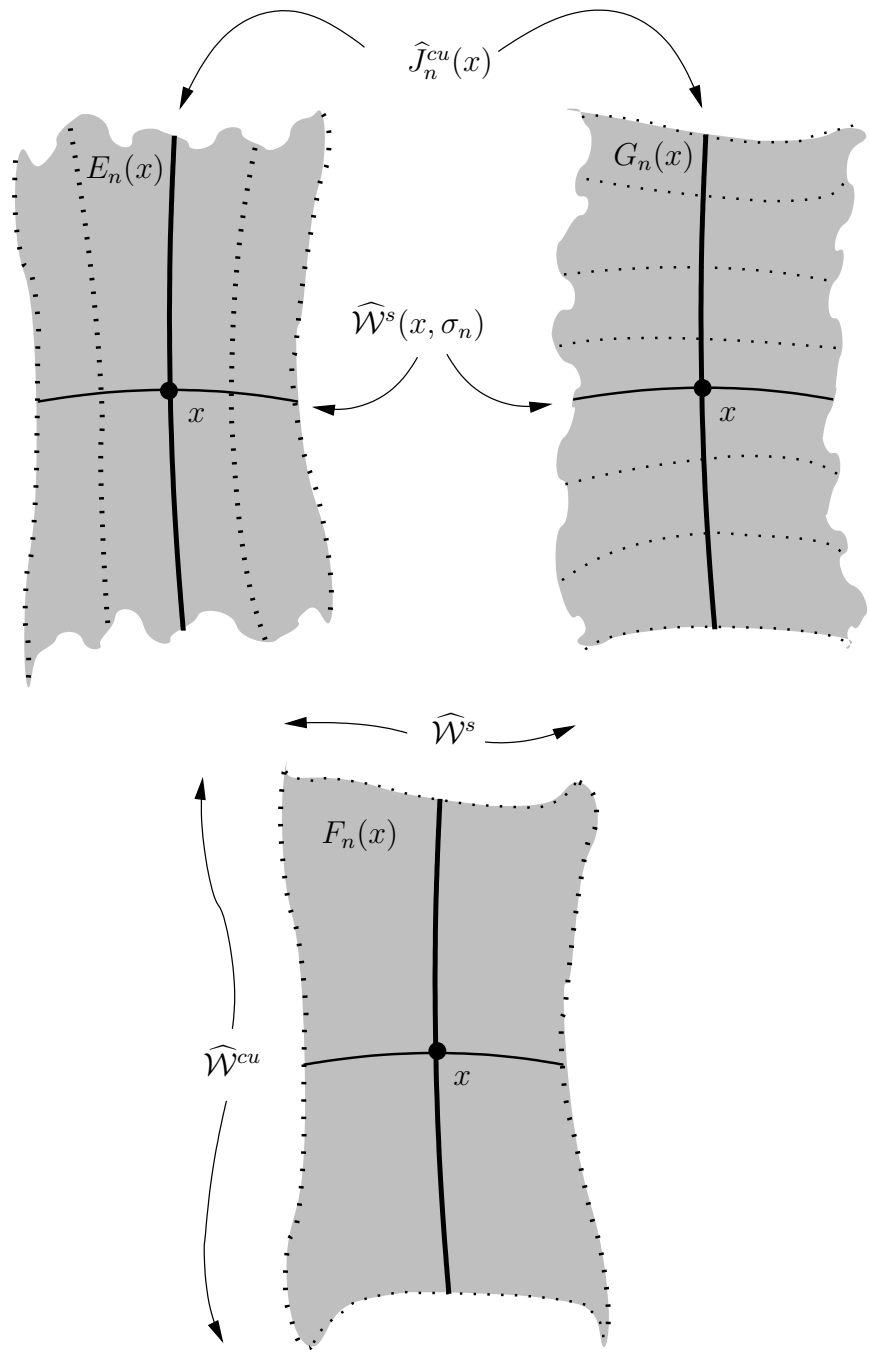} 
\end{center}
 \caption{Comparison between $E_n(x)$, $F_n(x)$ and $G_n(x)$.}
\end{figure}
Figure~\ref{f=efg} is a schematic illustration of the relationship between the sets
$E_n(x)$, $F_n(x)$ and $G_n(x)$. All three sets contain $\hJ^{cu}_n(x)$ and
$\W^s(x,\sigma_n)$. The set $E_n(x)$ fibers over $\W^s(x,\sigma_n)$ with
fibers of the form $\hJ^{cu}_n(\cdot)$. The set $G_n(x)$ fibers over 
$\hJ^{cu}_n(x)$ with fibers of the form $\W^s(\cdot,\sigma_n)$. The foliation
product $F_n(x)$ of $\hJ^{cu}_n(x)$ and $\W^s(x,\sigma_n)$ 
is, in some sense,
intermediate 
between $E_n(x)$ and $G_n(x)$.

We now prove these equivalences, following the outline described above.

First, recall that $B_n(x)$ is a round ball about $x$ of radius $\sigma_n$. 
The forward implication in the first equivalence  
is obvious from the definition of $B_n(x)$.
The backward implication follows from this definition and the fact that 
the ratio $\sigma_{n+1}/\sigma_{n} = \sigma(p_n)$
of successive radii is less than $1$, and is
bounded away from both $0$ and $1$ independently 
of $n$. From this we also see that $B_n(x)$ is regular.

The set $C_n(x)$ fibers over $D^{cs}_n(x)$, with fiber
$\W^u(x',\sigma_n)$ over $x'\in D^{cs}_n(x)$. 
The sequence $D^{cs}_n(x)$ is comparable to the 
sequence of disks $\hW^{cs}(x, \sigma_n)$, by continuity
of the foliation $\hW^c$.  Continuity of the foliation $\W^s$ then implies
that $C_n(x)$ is comparable to $B_n(x)$.

To prove the equivalence 
$$\lim_{n\to\infty} m(X:C_n(x)) = 1 \Longleftrightarrow
              \lim_{n\to\infty} m(X:D_n(x)) = 1,$$
we note that $C_n(x)$ and $D_n(x)$ both fiber over $D^{cs}_n(x)$,
with $\W^u$-fibers.  Since $X$ is essentially $\W^u$-saturated, 
Proposition~\ref{p=compmeas} implies that it suffices to show that the
fibers of $C_n(x)$ and $D_n(x)$ are both $c$-uniform.  The fibers
of of $C_n(x)$ are easily seen to be uniform, because they are
all comparable to balls in $\W^u$ of fixed radius $\sigma_n$.  
The fibers of $D_n(x)$ are the unstable juliennes $J^u_n(x')$,
for $x'\in D^{cs}_n(x)$.  Uniformity of these fibers follows
from Proposition~\ref{p=julmeas}.
 
We next prove:

\begin{lemma} The sequences $D_n(x)$ and $E_n(x)$ are comparable.
\end{lemma}
\proof  If $y\in D_n(x)$, then $y_n\in \W^u(q_n, \tau_n)$, for
some $q\in D^{cs}_n(x)$. The point $q$ is the unique intersection of
$\W^u(y)$ with $\hW^{cs}(x)$.  Similarly, if $y\in E_n(x)$, 
then $y_n\in \hW^u(q_n', \tau_n)$, for some $q'\in D^{cs}_n(x)$. 
The point $q'$ is the unique intersection of
$\hW^u(y)$ with $\hW^{cs}(x)$. 

Let $y\in D_{n}(x)$, let $q$ be the unique point of intersection 
of $\W^u(y)$ with $\hW^{cs}(x)$, and let $q'$ 
be the unique point of intersection 
of $\hW^u(y)$ with $\hW^{cs}(x)$.  To show that
$D_n(x)$ is comparable to $E_n(x)$, it suffices to show that
$$d(y_n,q_n') = O(\tau_n)$$
and 
$$d(q,q') = O(\sigma_n).$$
By Proposition~\ref{p=localfol}, part 1., $T\hW^u$ lies in the
$\eps$-cone about $E^u$; since $\W^u$ is everywhere tangent to 
$E^u$, it follows
that 
$$d(y_n, q_n') = O (d(y_n, q_n)) = O(\tau_n).$$
From this and the triangle inequality it follows that 
$$d(q_n, q_n') = O(\tau_n).$$
Since the leaves of both $\W^u$ and $\hW^u$ are uniformly contracted by
$f^{-1}$, we have that
$$d(q_j, q_j') \leq d(q_j, y_j) + d(y_j, q_j') = 
O(\hat\nu_{-j}^{-1}(p_n)\tau_n),$$
for $0 \leq j\leq n$.  This implies that 
$$d(y_n, q_n') = O(\tau_n).$$

Since $d(q, q') = O(\hat\nu_n\tau_n)$ and $\hat\nu\tau < \sigma$, we
obtain that $d(q, q') = O(\sigma_n)$.\endproof

We next show:
\begin{lemma} $E_n(x)$ is comparable to $F_n(x)$, and $F_n(x)$ is comparable to $G_n(x)$.  
\end{lemma}
\proof The comparability of $E_n(x)$ and
$F_n(x)$ is an immediate consequence of Lemma~\ref{l=holon}.

To see that $F_n(x)$ and
$G_n(x)$ are comparable, suppose that
$q'$ lies in the boundary of the fiber of $F_n(x)$ that lies in 
$\W^s(q)$
for some $q \in \hJ^{cu}_n(x)$. 
Then $q' \in \hJ^{cu}_n(x')$ for a point $x'$ that lies in
the boundary of $\W^s(x,\sigma_n)$. The diameters of 
$\hJ^{cu}_n(x)$ and $\hJ^{cu}_n(x')$
are both $O(\sigma_n)$, and $d(x,x') = \sigma_n$.
Hence, if $k$ is large enough, we will have
 $$
\sigma_{n+k} \leq d(q,q') \leq {\sigma_{n}}.
$$
Thus all points on the boundary of the fiber of $F_n(x)$  in $\W^s_{loc}(q)$
lie outside $\W^s(q,\sigma_{n+k})$ and inside $\W^s(q,\sigma_{n-k})$.
\endproof

We now know that $D_n(x), E_n(x), F_n(x)$ and $G_n(x)$ are all comparable.
As discussed above, to prove the fourth through sixth 
equivalences, it now suffices to show:
                                                            
\begin{lemma}\label{l=gnreg}
The sequence $G_n(x)$ is regular for each $x \in \W^s(p,1)$.
\end{lemma}
                                     
\proof 
The set
$$G_n(x) = \bigcup_{q\in \hJ^{cu}_n(x)} \W^s(q,\sigma_n)$$
fibers over $\hJ^{cu}_n(x)$, with $\W^s$-fibers 
$\W^s(q,\sigma_n)$.  Since $\W^s$ is absolutely continuous,
Proposition~\ref{p=unifreg} implies that regularity of $G_n(x)$ follows from 
regularity of the base sequence and fiber sequence.
Proposition~\ref{p=julmeas} implies that the sequence
$\hJ^{cu}_n(x)$ is regular in the induced measure $\hm_{cu}$.
As we remarked above, the ratio
$\sigma_{n+1}/\sigma_{n} = \sigma(p_n)$ is uniformly bounded below away
from $0$.  Consequently, the ratio
$$\frac{m_s(\W^s(q,\sigma_{n+1}))}{m_s(\W^s(q,\sigma_{n}))}$$
is bounded away $0$, uniformly in $x, q$, and $n$.
The regularity of $G_n(x)$ now follows from Proposition~\ref{p=unifreg}.
\endproof                                                                                                                         

To prove the final equivalence, we use the fact that
$G_n(x)$ fibers over $\hJ^{cu}_n(x)$ with $c$-uniform fibers 
and apply Proposition~\ref{p=compmeas2}.  Here we use the
fact that $X$ is $\W^s$-saturated.
This completes the proof of Proposition~\ref{p=sliceden}. \endproof

\end{document}